\documentclass[centertags,leqno]{article}
\usepackage{def}
\title{Reduction and lifting of special metacyclic covers}
\author{Stefan Wewers\\[1ex]
        University of Pennsylvania}
\begin{document}

%\selectlanguage{english}

\maketitle

%\selectlanguage{french}

%\begin{abstract}
%  On definit des rev\^etements speciaux comme des rev\^etements
%  metacyclic de la droite projective ayant un certain type de
%  r\'eduction mauvaise. Ces rev\^etement apparaissent naturellement dans
%  l'\'etude des rev\^etements galoisiens \'etales de la droite
%  projective moins trois points. Nos r\'esultats donnent une
%  classification compl\`ete des rev\^etements speciaux comme rel\`evement
%  canonique de certains donn\'ees en charact\'eristique $p$.
%\end{abstract}

%\selectlanguage{english}
      
\begin{abstract} 
  Special covers are metacyclic covers of the projective line, with
  Galois group $\ZZ/p\rtimes\ZZ/m$, which have a specific type of bad
  reduction to characteristic $p$. Such covers arise in the study of
  the arithmetic of Galois covers of $\PP^1$ with three branch points.
  Our results provide a classification of all special covers in terms of
  certain lifting data in characteristic $p$.
\end{abstract}

%------------------------------------------------------------------

\section*{Introduction}

In \cite{Raynaud98} Raynaud has given a criterion for good reduction
of Galois covers of the projective line which are ramified at three
points. The proof of this criterion depends on the analysis of Galois
covers with bad reduction, under certain conditions on the Galois
group. In one particular step of this analysis, Raynaud introduced the
notion of the {\em auxiliary cover}: to a $G$-Galois cover
$f:Y\to\PP^1$ with bad reduction to characteristic $p$ he associates
(under certain conditions) an $H$-Galois cover $f\aux:Y\aux\to\PP^1$,
which has, in some sense, the same type of bad reduction as $f$, but
whose Galois group $H$ is a certain solvable quotient of a subgroup of
$G$. For instance, if $p$ strictly divides the order of $G$, $H$ is a
metacyclic group, isomorphic to $\ZZ/p\rtimes\ZZ/m$. Thus, for many
purposes, the study of bad reduction of Galois covers can be reduced
to the study of covers with certain solvable (in the easiest case,
metacyclic) Galois groups. However, this reduction step is paid for by
the introduction of extra branch points (the branch locus of $f$ is a
subset of the branch locus of $f\aux$). In general, it is hard to
predict where these extra branch points occur. 

The present paper is concerned with a detailed study of metacyclic
covers of $\PP^1$, with Galois group isomorphic to
$\ZZ/p\rtimes\ZZ/m$, which arise as the auxiliary cover of $G$-covers
of $\PP^1$ with three branch points and prime-to-$p$ ramification. In
Section \ref{sec1}, we give a characterization of such covers which is
independent of the group $G$; this characterization gives rise to the
definition of {\em special metacyclic covers}. The results of Section
\ref{sec2} and Section \ref{sec3} give, in some sense, a complete
classification of all special covers and a thorough understanding of
their arithmetic. These results are applied in \cite{badred} to the
study of the arithmetic of three point Galois covers of $\PP^1$ with
bad reduction. In particular, upper and lower bounds on the
ramification index of $p$ in the field of moduli of such covers are
given, sharpening the results of \cite{Raynaud98}.

\vspace{1ex} Let us now give a more detailed outline of our results.
Fix a prime number $p$ and a finite, center free group $G$ such that
$p$ strictly divides the order of $G$. Let us also fix a field $K$ of
characteristic $0$ which is complete with respect to a discrete
valuation $v$. We assume that the residue field $k$ of $v$ is
algebraically closed and of characteristic $p$.  Suppose we are given
a $G$-Galois cover $f:Y\to\PP^1$ defined over $K$, branched at
$x_1=0$, $x_2=1$ and $x_3=\infty$, with prime-to-$p$ ramification. We
assume that $f$ has bad reduction at $v$, and let $\fb:\Yb\to\Xb$ be
the {\em stable reduction} of $f$ (for a precise definition, see
Section \ref{sec1.3}).

Let $f\aux:Y\aux\to\PP^1$ be the auxiliary cover associated to $f$.
This is a certain $H$-Galois cover, branched at $r\geq 3$ points
$x_1,\ldots,x_r$, where $x_1,x_2,x_3$ are the branch points of $f$ and
$H$ is the quotient of a subgroup of $G$, isomorphic to
$\ZZ/p\rtimes\ZZ/m$ for some integer $m>1$ dividing $p-1$. The cover
$f\aux$ has bad reduction, as well, and its stable reduction
$\fb\aux:\Yb\aux\to\Xb$ is is closely related to the stable reduction
of $f$ (e.g.\ $\fb$ and $\fb\aux$ have the same target $\Xb$). We
shall call an $H$-cover of $\PP^1$ {\em special} if it arises as the
auxiliary cover of a three point cover with prime-to-$p$ ramification.
(In Section \ref{sec1}, we will actually use a definition of specialty
that is independent of the existence of such a three point cover).

The natural question to ask now is: what do we know about the extra
branch points $x_4,\ldots,x_r$? After all, $f\aux$ is essentially
determined by the three point cover $f$, which is itself determined by
finite data. Therefore, the points $x_4,\ldots,x_r$ should be in some
kind of special position. Our first main result is:

\vspace{2ex}\noindent
{\bf Theorem A
  \em The (semistable) curve $\Xb$ is the union of $r+1$ distinct
  irreducible components $\Xb_0,\ldots,\Xb_r$ (each isomorphic to
  $\PP^1_k$) such that, for $i=1,\ldots,r$, $\xb_i\in\Xb_i$ and
  $\Xb_i$ meets $\Xb_0$ in a unique point $\tau_i$. In particular,
  $x_i\not\equiv x_j\mod{v}$ for $i\neq j$ (as points on $\PP^1$).
}\vspace{2ex}
    
Note that the isomorphism $\Xb_0\cong\PP^1_k$ is canonical and that
the point $\tau_i\in\PP^1_k$ is nothing but the specialization of
$x_i$, on the standard model of $\PP^1$. 

Theorem A essentially says that the stable reduction of a three point
Galois cover of the projective line is as simple as one can expect it
to be. This result turns out to be very useful to study the arithmetic
of such covers, see \cite{badred}.  Let us mention at this point that
the `three point condition' is essential for Theorem A to hold.

\vspace{1ex} The proof of Theorem A shows that the stable reduction
$\fb\aux:\Yb\aux\to\PP^1_k$ essentially corresponds to a certain pair
$(\Zb_0,\omega_0)$. Here $\Zb_0\to\PP^1_k$ is an $m$-cyclic cover
branched at the points $\tau_1,\ldots,\tau_r$ (as in Thm.\ A), and
$\omega_0$ is a regular differential form on $\Zb_0$. The differential
form $\omega_0$ verifies the following conditions: a.) $\omega_0$ is
logarithmic, i.e.\ of the form ${\rm d}u/u$, b.) $\omega_0$ is an
eigenvector under the action of $\ZZ/m$ and c.) the zeros of
$\omega_0$ are contained in the ramification locus of
$\Zb_0\to\PP^1_k$. Let us call the pair $(\Zb_0,\omega_0)$ a {\em
  special degeneration datum}. It turns out that Condition b.) and c.)
already determine the differential $\omega_0$ up to a constant factor.
Therefore, Condition a.) poses a restriction on the cover
$\Zb_0\to\PP^1_k$ and, in particular, on the branch points
$\tau_1,\ldots,\tau_r$. In concrete terms, this condition translates
into an explicit system of equations satisfied by the tuple
$(\tau_i)$.  One can show that this system has only a finite number of
solutions, and that therefore there exist only a finite number of
nonisomorphic special degeneration data, for a given prime $p$.

Our second main result is that every special degeneration datum arises
from a special $H$-cover. More precisely, we have:

\vspace{2ex}\noindent
{\bf Theorem B  
  \em Let $(\Zb_0,\omega_0)$ be a special degeneration datum, defined 
  over an algebraically closed field $k$ of characteristic $0$.  Let 
  $K_0$ be the fraction field of the ring $W(k)$ of Witt vectors over 
  $k$, and let $K:=K_0(\zeta_p^{(m)})/K_0$ be the unique extension of
  degree $(p-1)/m$ contained in $K_0(\zeta_p)$. Then there exists a
  special $H$-cover $f:Y\to\PP^1$, defined over $K$, whose reduction
  $\fb:\Yb\to\Xb$ corresponds to $(\Zb_0,\omega_0)$. The cover $f$ is
  unique up to the choice of the branch points
  $x_1,\ldots,x_r\in\PP^1(K_0)$ (which necessarily lift
  $\tau_1,\ldots,\tau_r\in\PP^1(k)$).
}\vspace{2ex}

This theorem gives us a concrete way to construct special $H$-covers
and, together with Theorem A, yields a complete classification of
special covers (up to solving the equations satisfied by the points
$\tau_i$).  Again, this result has nice applications to the arithmetic
of three point covers, see \cite{badred}.

Theorem B is a result on lifting covers of curves from characteristic
$p$ to characteristic $0$. In this sense, it is similar to the main
result of \cite{Yannick}. But even though \cite{Yannick} has
influenced our work a lot (this is most obvious in Section
\ref{sec2}), the actual proof of Theorem B uses a somewhat different
approach. The main differences are the following. First, the results
of \cite{Yannick} are formulated in terms of automorphisms of the
$p$-adic open disk (and are therefore something local), whereas we are
dealing with projective curves. Second, for the applications in
\cite{badred} it is important to prove uniqueness of lifting and to
have a good control over the field of definition. It seems to the
author that the methods of \cite{Yannick} are not sufficient to prove
such a result. Briefly, what helps us out is the extra $\ZZ/m$-action,
which rigidifies the lifting process.

\vspace{2ex}
The author gratefully acknowledges financial support from the {\em
  Deutsche Forschungsgemeinschaft}.

%-------------------------------------------------------------------

\section{Special metacyclic covers}  \label{sec1}

A metacyclic cover of the type we consider in this paper is the
composition of a ramified $m$-cyclic cover $Z\to\PP^1$ and an \'etale
$p$-cyclic cover $Y\to Z$. It therefore corresponds to an element of a
certain isotypical component of the $p$-torsion of the Jacobian $J_Z$.
We say that a metacyclic cover is of {\em multiplicative type} if this
isotypical component of $J_Z[p]$ has \'etale reduction. This forces
the cover $Y\to Z$ to have bad (more specifically, multiplicative)
reduction, and allows us to compute the field of moduli (Section
\ref{sec1.2}). In Section \ref{sec1.3}, we define {\em special covers}
as metacyclic covers $f:Y\to\PP^1$ of multiplicative type whose stable
reduction satisfies certain (rather restrictive) conditions. In
Section \ref{sec1.4}, we show that the auxiliary cover of a Galois
cover of $\PP^1$ branched at three points is special, provided that
$p$ strictly divides the order of the Galois group.

\subsection{Metacyclic covers of multiplicative type} \label{sec1.1}

Let us fix the following objects:
\begin{itemize}
\item 
  A complete discrete valuation ring $R_0$, with fraction field
  $K_0$ of characteristic $0$ and residue field $k$ of characteristic
  $p>0$.  We will assume that $k$ is algebraically closed. Let $\Kb$
  be an algebraic closure of $K_0$ and $\zeta_p\in\bmu_p(\Kb)$ a fixed
  $p$th root of unity. We write $v_{R_0}$ for the valuation of $K_0$
  corresponding to $R_0$. We normalize $v_{R_0}$ such that
  $v_{R_0}(p)=1$. 
\item
  An integer $m>0$ dividing $p-1$, and a character
  $\chi:\ZZ/m\to\bmu_m(K_0)$ of order $m$. We write
  $\chib:\ZZ/m\to\FF_p^\times$ for the reduction of $\chi$ mod $v_{R_0}$.
  We define the group $H:=\ZZ/p\rtimes\ZZ/m$ such that $\ZZ/m$ acts on
  $\ZZ/p$ via $\chib$.
\item
  An integer $r\geq 3$ and an $r$-tuple $(a_1,\ldots,a_r)$ of
  integers such that $0<a_i<m$ and $\sum_i a_i \equiv
  0\mod{m}$. 
\item 
  An $r$-tuple $(x_1,\ldots,x_r)$ of pairwise distinct
  $K_0$-rational points on the projective line $\PP^1$.  
\end{itemize}

The inverse image of an element $a\in(\ZZ/m)^\times$ under the natural
map $H\to\ZZ/m$ is a conjugacy class of $H$, containing elements of
order $m/(m,a)$. We denote this conjugacy class by $C_a$.

\begin{defn} \label{def1.1}
  A {\em metacyclic cover of type $(x_i;a_i)$}
  is an $H$-cover
  \[
       f:Y\lpfeil{H}\PP^1_{\Kb}
  \]
  defined over $\Kb$, branched in $x_1,\ldots,x_r$, such that the
  canonical generator of inertia above $x_i$ (with respect to the
  character $\chi$) is an element of $C_{a_i}$, for $i=1,\ldots,r$
  (see e.g.\ \cite{SerreTopics}).
  The cover $f$ is called {\em of  multiplicative type} if
  \begin{equation} \label{eq1.0}
             \sum_{i=1}^r \;a_i \;=\; m.
  \end{equation}
\end{defn}

The meaning of the term `of multiplicative type' will become clear
later. For the moment, let us fix a metacyclic cover $f:Y\to\PP^1$ of
type $(x_i;a_i)$, not necessarily of multiplicative type. Let $Z$ be
the quotient of $Y$ by the normal subgroup $\ZZ/p\lhd H$. Thus,
$g:Z\to\PP^1_{\Kb}$ is an $m$-cyclic cover, branched in
$x_1,\ldots,x_r$ such that $a_i$ is the canonical generator of inertia
above $x_i$.  We shall call $g:Z\to\PP^1$ the $m$-cyclic cover of type
$(x_i;a_i)$. In concrete terms, $Z$ is the normalization of the plane
curve with equation
\begin{equation} \label{eq1.1}
        z^m \;=\; \prod_{i=1}^r\; (x-x_i)^{a_i}
\end{equation}
(provided $x_i\not=\infty$), and $a\in\ZZ/m$ acts
on $Z$ as the automorphism $\phi_a:Z\iso Z$ such
that $\phi_a^*z=\chi(a)\,z$. 

The $p$-cyclic cover $Y\to Z$ is \'etale. Hence it corresponds to a
nontrivial class $\theta$ in $H_{\rm \acute{e}t}^1(Z,\ZZ/p)_{\chi}$
(for any $\FF_p[\ZZ/m]$-module $M$, we denote by $M_{\chi}$ the
$\chi$-eigenspace). Let $J_Z={\rm Pic}^0(Z)$ denote the Jacobian of
$Z$, and $J_Z[p]$ its group of $\Kb$-rational points of order $p$. 
Kummer theory gives us a canonical isomorphism
\begin{equation} \label{eq1.3} 
   H_{\rm \acute{e}t}^1(Z,\ZZ/p)_{\chi} \;\cong\; J_Z[p]_\chi(-1)
   \;:=\; \Hom_{\FF_p}(\bmu_p(\Kb),J_Z[p]_\chi),
\end{equation}
see \cite{Milne}. Using our choice of a $p$th root of unity $\zeta_p$,
we can identify $J_Z[p]_\chi(-1)$ with $J_Z[p]_\chi$. Let $L$ be the
line bundle on $Z$ of order $p$ corresponding to the class $\theta$
under the isomorphism \eqref{eq1.3}. Let $D$ be a divisor on $Z$ such
that $L\cong\O_Z(D)$. Then $p\cdot D$ is the divisor of some rational
function $u$ on $Z$. By the definition of the isomorphism
\eqref{eq1.3}, the cover $Y\to Z$ is (birationally) given by the
equation
\begin{equation} \label{eq1.4}
    y^p \;=\; u,
\end{equation}
and $b\in\ZZ/p$ acts on $Y$ as the automorphism $\psi_b:Y\iso Y$ such
that $\psi_b^*y=\zeta_p^b\,y$. One can easily check that the action of
$\ZZ/m\subset H$ on $Y$ is obtained by extending the automorphisms
$\phi_a$, $a\in\ZZ/m$, from $Z$ to $Y$ such that
$\phi_a^*y:=\zeta_p^{\chib(a)}y$.

\subsection{The field of moduli} \label{sec1.2}

Let $K\inn$ be the {\em field of moduli} of the $H$-cover
$f:Y\to\PP^1$, relative to the extension $\Kb/K_0$. Since the group
$H$ has trivial center, this means that $K\inn/K_0$ is the smallest
field extension such that the $H$-cover $f$ descents to an $H$-cover
$f_{K\inn}:Y_{K\inn}\to\PP^1_{K\inn}$ over $K\inn$.  Moreover, the
extension $K\inn/K_0$ is finite, and the model $f_{K\inn}$ of $f$ is
unique up to $K\inn$-linear isomorphism. We let
$\Gamma\inn:=\Gal(\Kb/K\inn)$.

The field of moduli of the $m$-cyclic cover $g:Z\to\PP^1$ is just
$K_0$. Although $g$ has no unique $K_0$-model (because $\ZZ/m$ is
abelian), there is a canonical model $g_{K_0}:Z_{K_0}\to\PP^1_{K_0}$,
given by the equation \eqref{eq1.1}. (Since $K_0$ is strictly
henselian, $g_{K_0}$ is characterized by the fact that it is
unramified at the generic point of $\PP^1_k$. Note, however, that
the $m$-cyclic quotient $Y_{K\inn}/(\ZZ/p)\to\PP^1_{K\inn}$ of
$f_{K\inn}$ will not be $K\inn$-isomorphic to $g_{K_0}\otimes K\inn$,
in general.) The choice of the model $g_{K_0}$ determines an action of
$\Gamma_{K_0}:=\Gal(\Kb/K_0)$ on the $\FF_p[\ZZ/m]$-module
$H^1(Y,\ZZ/p)$, and hence an action on the $\chi$-eigenspace
$H^1(Z,\FF_p)_\chi$. We can describe the subgroup
$\Gamma\inn\subset\Gamma_{K_0}$, and therefore the field extension
$K\inn/K_0$, in terms of this action, as follows.

\begin{prop} \label{prop1.1}
  Let $\theta\in H^1(Z,\FF_p)_\chi$ be the class corresponding to the 
  \'etale $p$-cyclic cover $Y\to Z$. Then
  \[
       \Gamma\inn \;=\; \{\;\sigma\in\Gamma_{K_0} \mid
          \op{\sigma}{\theta}=\chib(a)\,\theta,\;
            a\in\ZZ/m\;\}.
  \]
\end{prop}

\proof Let $\sigma\in\Gamma_{K_0}$. By definition, $\sigma$ is an
element of $\Gamma\inn$ if and only if the conjugate cover
$\op{\sigma}{f}:\op{\sigma}{Y}\to\PP^1$ is isomorphic to $f$.  This
means that there exists a $\Kb$-linear isomorphism
$\phi:\op{\sigma}{Y}\iso Y$ which is equivariant with respect to the
$H$-action, such that $f\circ\phi=\op{\sigma}{f}$. We may (and will)
identify $\op{\sigma}{Z}$ with $Z$ (using the model $Z_{K_0}$). So, if
it exists, the isomorphism $\phi$ restricts to a $\Kb$-linear
automorphism of $Z$ which commutes with the $\ZZ/m$-action and the map
$g:Z\to\PP^1$; therefore, $\phi|_Z=\phi_a$, for some $a\in\ZZ/m$. In
other words, $\sigma\in\Gamma\inn$ if and only if
$\op{\sigma}{\theta}=\phi_a^*\theta=\chib(a)\theta$, for some
$a\in\ZZ/m$. This is what we wanted to prove.  \Endproof

Of course, the choice of the model $g_{K_0}$ also determines an action
of $\Gamma_{K_0}$ on $J_Z[p]_\chi$ and $J_Z[p]_\chi(-1)$, such that
\eqref{eq1.3} becomes an isomorphism of 
$\FF_p[\Gamma_{K_0}]$-modules. One way to study this action is to
regard $J_Z[p]_\chi$ as the group of $\Kb$-rational points of the {\em
group scheme} $J_{Z_{K_0}}[p]_\chi$ over $K_0$. A crucial fact we will
use in this paper is the following.

\begin{prop} \label{prop1.2}
  Suppose that $f$ is of  multiplicative type, i.e.\ Equation
  \eqref{eq1.0} holds. Then
  \begin{enumerate}
  \item 
    As a group scheme over $K_0$, $J_{Z_{K_0}}[p]_\chi$ is isomorphic to
    $(\ZZ/p)^{r-2}$. Therefore, the action of $\Gamma_{K_0}$ on
    $J_Z[p]_\chi$ is trivial.
  \item
    We have $K\inn=K_0(\zeta_p^{(m)}):=K_0(\zeta_p)^{\ZZ/m}$ (where 
    $a\in\ZZ/m$ acts on $K_0(\zeta_p)$ as
    $\zeta_p\mapsto\zeta_p^{\chib(a)}$).
  \end{enumerate}
\end{prop}

\proof Part (i) is proved in \cite{IreneAux}, Proposition
1.6. We will outline the proof in a special case in Section
\ref{sec1.2b} below.  Part (ii) follows from the first and
Proposition \ref{prop1.1}.  \Endproof

\begin{cor} \label{cor1.1}
  If $f$ is of multiplicative type, then it has bad reduction.
\end{cor}

\proof See \cite{IreneAux}, Proposition 2.9 (and also the next
subsection).  \Endproof

\subsection{If the branch points do not coalesce} \label{sec1.2b}

Let $f:Y\to\PP^1$ be a metacyclic cover of multiplicative type
$(x_i;a_i)$. In this subsection, we assume that $x_i\not\equiv
x_j\mod{v_{R_0}}$, for $i\neq j$. In other words, we assume that the
branch points $x_i$ specialize to pairwise distinct $k$-rational
points $\xb_i\in\PP^1_k$. In this special case, we would like to
explain some of the ingredients of the proof of Proposition
\ref{prop1.2} (i) and Corollary \ref{cor1.1}, because they will be
useful later on.

As before, we let $g_{K_0}:Z_{K_0}\to\PP^1_{K_0}$ be the canonical
model over $K_0$ of the $m$-cyclic quotient $g:Z\to\PP^1$ of $f$.  By
assumption, $g_{K_0}$ extends to a tamely ramified $m$-cyclic cover
$g_{R_0}:Z_{R_0}\to\PP^1_{R_0}$, with $Z_{R_0}$ smooth over $R_0$. The
special fiber $\gb:\Zb\to\PP^1_k$ of $g_{R_0}$ is the $m$-cyclic cover
of type $(\xb_i;a_i)$. By this we mean that $a_i\in\ZZ/m$ is the
canonical generator of inertia above $\xb_i$, with respect to the
character $\chib:\ZZ/m\to k^\times$. The theory of N\'eron models
defines a surjective specialization map
\begin{equation} \label{eq1.4b}
   J_Z[p]_{\chib}  \To  J_{\Zb}[p]_{\chib},
\end{equation}
which we regard as a morphism of $\FF_p$-modules.  Note that
\eqref{eq1.4b} is an isomorphism if and only if the group scheme
$J_{Z_{K_0}}[p]_{\chib}$ extends to a finite \'etale group scheme over
$R_0$. Therefore, Proposition \ref{prop1.2} (i) (in the special case
we consider here) follows from the assertion that \eqref{eq1.4b} is an
isomorphism between $\FF_p$-modules of dimension $r-2$.  Using the
identification \eqref{eq1.3} and the Riemann--Hurwitz formula, one
shows that $\dim_{\FF_p}J_Z[p]_{\chib}=r-2$. Thus, Proposition
\ref{prop1.2}, in the case where $x_i\not\equiv x_j\mod(v_{R_0})$,
follows from Lemma \ref{lem1.0} (iii) below.

Let $\C:H^0(\Zb,\Omega^1_{\Zb/k})\to H^0(\Zb,\Omega^1_{\Zb/k})$ be the
{\em Cartier operator}. By \cite{SerreMex}, we have a canonical
isomorphism
\begin{equation} \label{eq1.4c}
  J_{\Zb}[p] \;\cong\; 
     H^0(\Zb,\Omega^1_{\Zb/k})^{\cal \textstyle C} \;:=\;
     \{\;\omega\in H^0(\Zb,\Omega^1_{\Zb/k}) \;\mid\; 
                    \C(\omega)=\omega\;\}.
\end{equation}
Here $J_{\Zb}[p]$ means the $\FF_p$-module of {\em physical}
$p$-torsion points. Since $\C$ is $p^{-1}$-linear and $m|p-1$, $\C$
restricts to an endomorphism of the $\chib$-eigenspace
$H^0(\Zb,\Omega^1_{\Zb/k})_{\chib}$.

\begin{lem} \label{lem1.0}
  Let $\xb_1,\ldots,\xb_r$ be pairwise distinct $k$-rational points on
  $\PP^1_k$, and let $\Zb\to\PP^1_k$ be the $m$-cyclic cover of type
  $(\xb_i;a_i)$ (with respect to $\chib$). Assume that $\sum_i
  a_i=m$. Then:
\begin{enumerate}
\item
   $\dim_k H^0(\Zb,\Omega^1)_\chi = r-2$,
\item
   $\C$ is a bijection on $H^0(\Zb,\Omega^1)_\chi$, and
\item
   $\dim_{\FF_p}J_{\Zb}[p]_{\chib}=r-2$.
\end{enumerate}
\end{lem}

\proof By Serre duality, $H^0(\Zb,\Omega^1)_\chi$ is dual to
$H^1(\Zb,\O_{\Zb})_{\chi^{-1}}$.  Via this duality, $\C$ corresponds
to the Frobenius operator $F$ on $H^1$. Therefore, Part (i) of the
lemma follows from \cite{IreneComp}, Lemma 1.3.4 (note that $m-a_i$ is
the canonical generator of inertia above $\xb_i$, with respect to the
character $\chib^{-1}$). Also, to prove Part (ii), it suffices to show
that $F$ is a bijection on $H^1(\Zb,\O_{\Zb})_{\chib^{-1}}$. It
follows from Theorem 2.2.1 of \cite{IreneComp} that this is true at
least for a generic choice of the branch points $\xb_i$ (more
precisely, for $(\xb_i)$ contained in some open and dense subset of
$(\PP^1_k)^r-\Delta_r$). Furthermore, Part (iii) of the lemma follows
from Equation \eqref{eq1.4c} and Part (i) and (ii) (use \cite{Milne},
Lemma III.4.13). This completes the proof of the lemma, in the case of
generic branch points $\xb_i$.

The general case follows from the generic case by a purity
argument. Here is a short outline. Let $\Z\to\PP^1_S$ be the
`universal family over $\ZZ_p$' of all $m$-cyclic covers of $\PP^1$ of
type $(a_i)$. The crucial fact here is that the base $S$ of this
family is regular. Let $\G:=J_{\Z}[p]_{\chib}$ be the
$\chib$-isotypical part of the $p$-torsion of the Jacobian of $\Z$;
this is a finite flat group scheme over $S$ (of degree $p^{r-2}$). The
`generic case' of Lemma \ref{lem1.0} (ii) shows that
$\G\otimes_{\ZZ_p}\FF_p$ is generically \'etale over
$S\otimes_{\ZZ_p}\FF_p$. Therefore, $\G$ is \'etale over $S$ outside a
subspace of codimension $\geq 2$. Purity of Branch Locus shows that
$\G/S$ is actually \'etale everywhere. This completes the proof of the
lemma in general.  \Endproof

The purity argument used in the proof of Lemma \ref{lem1.0} is also
the main ingredient to prove Proposition \ref{prop1.2} in full
generality, i.e.\ without the condition $x_i\not\equiv x_j$. See
\cite{IreneAux} for details.

Let us sketch the proof of Corollary \ref{cor1.1} in the case that
$x_i\not\equiv x_j\mod{v_{R_0}}$. The \'etale $p$-cyclic cover $Y\to
Z$, regarded as a $\bmu_p$-torsor, corresponds to a closed embedding
$\ZZ/p\inj J_Z$ of group schemes over $\Kb$. If $f$ is of
multiplicative type, then Proposition \ref{prop1.2} says that this
embedding extends to an embedding $\ZZ/p\inj J_{Z_{R_0}}$ of group
schemes over $R_0$. By \cite{Raynaudfest}, \S 3, this means that $Y\to
Z$ has {\em multiplicative reduction}, i.e.\ $Y\to Z$ extends to a
$\bmu_p$-torsor $Y'_R\to Z_R$ over $R:=R_0[\zeta_p]$, whose
restriction to the special fiber is inseparable. In particular, $f$
has bad reduction.

The proof in the general case is similar, using e.g.\ \cite{Saidi00}
instead of \cite{Raynaudfest} to analyze the reduction of $Y\to Z$ in
case $Z$ extends to a semistable $R$-curve $Z_R$. See also Section
\ref{sec2.1} and \ref{sec2.2}.

\subsection{Special covers}   \label{sec1.3}

For the rest of this section, we will assume that $f:Y\to\PP^1$ is a
metacyclic cover of multiplicative type $(x_i;a_i)$. We do {\em not}
assume that $x_i\not\equiv x_j\mod{v_{R_0}}$. By Corollary
\ref{cor1.1}, $f$ has bad reduction. Let $K\inn/K_0$ be the field of
moduli of $f$ and $f_{K\inn}:Y_{K\inn}\to\PP^1_{K\inn}$ the model of
$f$ over $K\inn$.  Furthermore, let $K/K\inn$ be a sufficiently large
finite extension over which $f_{K\inn}$ has {\em stable reduction},
see \cite{badred}.  Recall that this means the following. First, the
ramification points $y_1,\ldots,y_s\in Y$ of the $H$-cover $f$ are
$K$-rational points on $Y_{K\inn}$. Second, after the base change
$Y_K:=Y_{K\inn}\otimes_{K\inn} K$, the pointed curve $(Y_K;y_i)$
extends to a stably pointed curve $(Y_R;y_{R,i})$ over the ring of
integers $R$ of $K$. The quotient map $f_R:Y_R\to X_R:=Y_R/H$ is
called the {\em stable model} of $f$ over $R$. Its special fiber
$\fb:\Yb\to\Xb$ is called the {\em stable reduction} of the $H$-cover
$f$. Note that $\fb$ is a finite map between pointed semistable
$k$-curves. Here the distinguished points on $\Yb$ are the
specializations $\yb_1,\ldots,\yb_s$ of the ramification points, and
the distinguished points on $\Xb$ are the specializations
$\xb_1,\ldots,\xb_r$ of the branch points of $f$.  However,
$(\Xb;\xb_i)$ is not a stably pointed curve, in general.

We call a component $\Xb_v$ of $\Xb$ {\em separable} if the
restriction of $\fb$ to one (and therefore to all) components of $\Yb$
above $\Xb$ is a separable morphism. Note that $\Xb$ has at least one
component that is not separable, by the definition of bad reduction.
Let us denote by $\Xb_i$ the component of $\Xb$ containing $\xb_i$.

For convenience, we make the following assumption, which will be part
of the definition of a special cover (Definition \ref{specialdef}
below).

\begin{ass} \label{specialass}
  \begin{enumerate}
  \item 
    A component $\Xb_v$ of $\Xb$ is separable if and only if it is a
    tail, i.e.\ if it is connected to the rest of $\Xb$ in a single point.
  \item
    The components $\Xb_1,\ldots,\Xb_r$ are pairwise distinct, and
    they are precisely the tails of $\Xb$. 
  \end{enumerate}
\end{ass}

\begin{rem} \label{rem1.1}
  It is easy to see that the components $\Xb_i$ are separable.  If the
  branch points of $f$ do not coalesce, i.e.\ $x_i\not\equiv
  x_j\mod{v}$ for $i\not=j$, then part (i) of Assumption
  \ref{specialass} holds automatically, by \cite{Raynaud98}, Lemme
  3.1.2. In this case, Assumption \ref{specialass} (ii) can be phrased
  (in the terminology of \cite{Raynaud98}) as: `there are no new
  tails'. One should keep in mind that this is a rather strong
  condition. 
\end{rem}

Fix $i\in\{1,\ldots,r\}$. The restriction of $\fb$ to the tail $\Xb_i$
is a (possibly disconnected) $H$-Galois cover $\fb_i:\Yb_i\to\Xb_i$.
The cover $\fb_i$ is ramified at exactly two points, namely at $\xb_i$
and at the unique intersection point of $\Xb_i$ with the rest of $\Xb$
(let us call this point $\tau_i$). The ramification above $\xb_i$ is
tame and, since $\xb_i$ is the specialization of the branch point
$x_i$ of $f$, the canonical generator of inertia above $\xb_i$ is in
the conjugacy class $C_{a_i}$.  On the other hand, the ramification at
$\tau_i$ is wild, of order $pm_i$, where $m_i:=m/(a_i,m)$.  The cover
$\fb_i:\Yb_i\to\Xb_i\cong\PP^1_k$ is of `Katz-Gabber-type'; see e.g.\ 
\cite{GilleLuminy}, \S 2, where such covers are called
`$m_i$-special'.  We set $\sigma_i:=h_i/m_i$, where $h_i$ is the
conductor of the $p$-part of inertia. Recall that $\sigma_i$ is the
jump in the inertia group filtration of the cover
$\fb_i:\Yb_i\to\Xb_i$, with respect to the {\em upper numbering} (see
\cite{SerreCL}).

\begin{prop} \label{prop1.3}\ 
    We have $\sigma_i=\nu_i+a_i/m$, for integers $\nu_i\geq 0$ such
    that 
    \begin{equation} \label{eq1.5}
          \sum_{i=1}^r \nu_i \;=\; r-3.
    \end{equation}
\end{prop}

\proof 
An easy argument shows that $h_i\equiv a_i/(m,a_i)\mod{m_i}$, see
\cite{IreneAux}. Therefore, $\sigma_i=\nu_i+a_i/m$, for integers
$\nu_i\geq 0$. On the other hand, we have Raynaud's vanishing cycle
formula, \cite{Raynaud98}, \S 3.4.2 (5). In our situation, it becomes
\begin{equation} \label{eq1.6}
    \sum_{i=1}^r\;  (1-\sigma_i)  \;=\;  2.
\end{equation}
Now the proposition follows from \eqref{eq1.6} and \eqref{eq1.0}.
\Endproof

Here comes the central definition of this paper.

\begin{defn} \label{specialdef}
  Let $f:Y\to\PP^1$ be a metacyclic cover of multiplicative 
  type $(x_i;a_i)$, with (bad) reduction $\fb:\Yb\to\Xb$. Suppose
  Assumption \ref{specialass} holds, and let $\sigma_i=\nu_i+a_i/m$ be
  the invariant attached to the wild ramification above the tail
  $\Xb_i$, as above. We say that $f$ is {\em special} if there are
  exactly three indices $i_1,i_2,i_3$ with $\nu_{i_j}=0$ (therefore,
  $\nu_i=1$ for $i\not\in\{i_1,i_2,i_3\}$, by Proposition
  \ref{prop1.3}).
\end{defn}

Thus, after reordering the branch points $x_i$, we may assume that   
$\nu_1,\nu_2,\nu_3=0$ and $\nu_i=1$ for $i=4,\ldots,r$.

\subsection{The auxiliary cover of a three point cover is special} 
\label{sec1.4}

In this subsection we discuss briefly the example that motivated our
study of special covers. For more details, and for arithmetic
applications of the results obtained in the present paper, see
\cite{badred}.

Let $G$ be a finite, center free group such that $p$ strictly divides
the order of $G$. Let $h:W\to\PP^1$ be a $G$-cover over $\Kb$,
ramified in $0$, $1$ and $\infty$, with ramification indices of order
prime-to-$p$.  Assume that the cover $h$ has bad reduction. Let
$h_R:W_R\to X_R$ be its stable model and $\hb:\bar{W}\to\Xb$ the
special fiber of $h_R$. By \cite{Raynaud98}, \S 3.1, the map $\hb$ is
separable exactly over the tails of $\Xb$. Let $\Xb_1$, $\Xb_2$,
$\Xb_3$ be the {\em primitive tails}, i.e.\ those which contain the
specialization of the branch points $0$, $1$, $\infty$, respectively.
Let $\Xb_4,\ldots,\Xb_r$ be the new tails. As above, we associate to
each tail $\Xb_i$ the invariant $\sigma_i=h_i/m_i$. By
\cite{Raynaud98}, Proposition 3.3.5, we have
\begin{equation} \label{eq1.7}
      \sigma_i >1, \qquad \text{for}\;\; i=4,\ldots,r.
\end{equation}
For $i=4,\ldots,r$, we choose a smooth point $\xb_i$ on $\Xb_i$, and
lift it to a point $x_i$ of $\PP^1_K$. Set $x_1:=0$, $x_2:=1$,
$x_3:=\infty$.
  
Now let $f:Y\to\PP^1$ be the {\em (small) auxiliary cover} associated
to $h$, see \cite{badred}. (Note that what is called the auxiliary
cover in \cite{Raynaud98} is called the {\em big} auxiliary cover in
\cite{badred}. The (small) auxiliary cover is defined as its quotient
by the prime-to-$p$ centralizer.) The cover $f$ is metacyclic of type
$(x_i;a_i)$, with $m$ and $(a_1,\ldots,a_r)$ as at the beginning of
Section \ref{sec1.1}. Furthermore, $f$ has bad reduction, and the
special fiber of its stable model is of the form $\fb:\Yb\to\Xb$
(where $\Xb$ is the target curve of $\hb$). Finally, the tails $\Xb_i$
of $\Xb$ are exactly the separable components (for the map $\fb$), and
the induced maps $\fb_i:\Yb_i\to\Xb_i$ over the tails are (possibly
disconnected) $H$-Galois covers of Katz-Gabber type, tamely ramified
at $\xb_i$ and wildly ramified at $\tau_i$, with invariant $\sigma_i$.
  
We claim that $f$ is special. Indeed, Assumptions \ref{specialass} (i)
and (ii) hold by construction (except maybe for $f$ being of
multiplicative type, i.e.\ $\sum_i a_i=m$). As in the proof of
Proposition \ref{prop1.3}, we show that $\sigma_i=\nu_i+a_i/m$, with
integers $\nu_i\geq 0$. Using again the vanishing cycle formula
\eqref{eq1.6}, we get
\begin{equation} \label{eq1.8}
      \sum_i \;\nu_i \;=\; r-2-\frac{\sum_i a_i}{m}.
\end{equation}
But $\nu_i>0$ for $i=4,\ldots,r$, by \eqref{eq1.7}, so the left hand
side of \eqref{eq1.8} is $\geq r-3$. We conclude that $\sum_ia_i=m$,
$\nu_i=0$ for $i=1,2,3$ and $\nu_i=1$ for $i=4,\ldots,r$. In other
words, $f$ is special.

\vspace{1ex}
It is very likely that every special $H$-cover $f:Y\to\PP^1$ arises as
the auxiliary cover of a $G$-cover $h:W\to\PP^1$ with three branch
points, for {\em some} group $G$. At the moment, we can prove this
claim only modulo the following hypothesis.

\begin{hyp} \label{specialhyp}
  Let $p$ be a prime and $m,h>1$ integers such that $m|p-1$, $(h,p)=1$
  and $1<\sigma:=h/m<2$. Then there exists a finite group $G$ and an
  \'etale $G$-cover $f:W\to\AA^1_k$ of the affine line, defined over
  $k=\bar{\FF}_p$, such that the inertia group at infinity is
  isomorphic to $\ZZ/p\rtimes\ZZ/m$ and the jump in the inertia group
  filtration is equal to $\sigma=h/m$. 
\end{hyp}

The hypothesis implies the claim, essentially because we can `reverse'
the auxiliary cover construction, using formal patching. This fact
establishes a strong link between the existence of \'etale Galois
covers of the affine line in characteristic $p$ with prescribed
ramification at infinity (i.e.\ certain generalized forms of
Abhyankar's Conjecture) and the reduction of Galois covers from
characteristic $0$ to characteristic $p$. For results in this
direction, see \cite{RRR}. These results show that we cannot
expect Hypothesis \ref{specialhyp} to hold for a given group $G$
satisfying the obvious restrictions (quasi-$p$, contains $H$).

%---------------------------------------------------------------------

\section{The structure of special covers} \label{sec2}

The main goal of this section is to prove the following result, which
essentially states that (the reduction of) special covers is as
simple as we can expect it to be.

\begin{thm} \label{thm1}
  Let $f:Y\to\PP^1$ be a special metacyclic cover of type $(x_i;a_i)$,
  with reduction $\fb:\Yb\to\Xb$. Then the curve $\Xb$ is the union of
  $r+1$ components $\Xb_0,\Xb_1,\ldots,\Xb_r$, such that $\Xb_0$ meets
  $\Xb_i$ in exactly one point, and $\Xb_i$ is the
  component to which $x_i$ specializes. In particular, $x_i\not\equiv
  x_j\mod{v_{R_0}}$ (as points on $\PP^1$) for $i\not=j$. 
\end{thm}

This was called Theorem A in the introduction. Here is a brief
outline of the proof.  Using results of \cite{Saidi00} (which
essentially reformulate earlier results of \cite{Raynaudfest},
\cite{GreenMatignon99} and \cite{Yannick} in a form suitable for us),
we analyze the stable reduction of the $p$-cyclic part $Y\to Z$ of a
given special cover $f:Y\to \PP^1$. We find that, in case Theorem
\ref{thm1} did not hold, one of the components of the special fiber
$\Yb$ of $Y$ would be an $\balpha_p$-torsor. Moreover, this
$\balpha_p$-torsor would have to satisfy certain strict numerical
conditions, coming from specialty of $f$ and compatibility with the
$\ZZ/m$-action. Then a direct calculation shows that such an
$\balpha_p$-torsor cannot exist, and Theorem \ref{thm1} follows.

In Sections \ref{sec2.1} and \ref{sec2.2} we analyze the stable
reduction of $f$ and translate its properties into a combinatorial
language similar to the notion of a {\em Hurwitz tree} developed in
\cite{Yannick}. Here we work under slightly more general assumptions
then is actually necessary for the proof of Theorem \ref{thm1}. This
will be useful in Section \ref{sec3}.

In Section \ref{sec2.3} we prove the key lemma which implies Theorem
\ref{thm1}. In Section \ref{sec2.4} we state a strengthening of
Theorem \ref{thm1}, which concerns the monodromy action on the special
fiber of a special cover.

\subsection{Admissible reduction of $Z\to\PP^1$} \label{sec2.1}

Let $f:Y\to\PP^1$ be a metacyclic cover of multiplicative type
$(x_i;a_i)$. Concerning the reduction of $f$, we use the notation
introduced in Section \ref{sec1.3}. In particular, $f_R:Y_R\to X_R$ is
the stable model of $f$, defined over a sufficiently large finite
extension $R/R_0$, with special fiber $\fb:\Yb\to\Xb$. The branch
point $x_i$ specializes to a point $\xb_i\in\Xb$, and $\Xb_i$ denotes
the component of $\Xb$ containing $\xb_i$. We make the following
assumption, which is slightly weaker then Assumption \ref{specialass}.
In particular, it holds if $f$ is special.

\begin{ass} \label{weakass}
\begin{enumerate}
\item
    The map $\fb:\Yb\to\Xb$ is separable exactly over the tails of $\Xb$.
  \item
    The components $\Xb_1,\ldots,\Xb_r$ are pairwise distinct tails.
  \end{enumerate}
\end{ass}

We define a tree $T$ associated to the semistable curve $\Xb$, as
follows.  The set of {\em vertices} $V$ of $T$ is defined as the set
of (irreducible) components of $\Xb$. For $v\in V$, we refer to the
corresponding component of $\Xb$ as $\Xb_v$. The set of {\em edges}
$E$ of $T$ is, by definition, the set of triples $e=(\xb_e,v,v')$, where
$\xb_e\in\Xb$ is a singular point and $v,v'\in V$ are vertices such
that $\xb_e\in\Xb_v\cap\Xb_{v'}$. The assignments $s(e):=v$ and
$t(e):=v'$ define the {\em source map} $s:E\to V$ and the {\em target
  map} $t:E\to V$. The edge $\bar{e}:=(\xb_e,v',v)\in E$ is called the
{\em opposite edge} of $e$. 

We define $B\subset V$ as the set of {\em leaves} of $T$, i.e.\ the
vertices corresponding to the tails of $\Xb$.  In other words, for
each $v\in B$ there is a unique edge $e_v$ such that $t(e_v)=v$. We
call the elements of $V':=V-B$ the {\em interior} vertices.  By
Assumption \ref{weakass} (ii), we may identify the set
$I:=\{1,\ldots,r\}$ with the corresponding subset of $B$. Note that
$B=I$ if $f$ is special.

Let us state the following easy lemma without proof.

\begin{lem} \label{lem2.1}
  There exists a unique map $E\to\{0,\ldots,m\}$, $e\mapsto a_e$,
  such that
  \begin{enumerate}
  \item 
    $a_e=a_i$, if $i:=t(e)\in I$, and $a_e=0$ if $v:=t(e)\in B-I$,
  \item 
    for all interior vertices $v$, we have $\sum_{s(e)=v}
    a_e \;=\; m$, and
  \item
    $a_{\bar{e}}+a_e=m$, for all edges $e\in E$.
  \end{enumerate}
  If $f$ is special, then $0<a_e<m$ for all $e\in E$.
\end{lem}

Let $Z_R:=Y_R/(\ZZ/p)$ be the quotient of $Y_R$ by the $p$-cyclic
normal subgroup of $H$. It is well known that $Z_R$ is again a
semistable curve over $R$, with generic fiber $Z_K$. Also, the action
of $\ZZ/m$ on $Z$ extends to $Z_R$, and the (tamely ramified)
$m$-cyclic cover $g_K:Z_K\to\PP_K$ extends to a finite morphism
$g_R:Z_R\to X_R$ of semistable $R$-curves. Since, by Assumption
\ref{weakass} (ii), the branch points $x_i$ of $g_K$ specialize to
pairwise distinct points of $\Xb$, the cover $g_R:Z_R\to X_R$ is an
$m$-cyclic {\em admissible} cover over $R$. In particular, its special
fiber $\gb:\Zb\to\Xb$ is an $m$-cyclic admissible cover over $k$. For
$v\in V$, we set $\Zb_v:=\gb^{-1}(\Xb_v)$.

\begin{prop} \label{prop2.1}
  For each vertex $v\in V$, the restriction $\gb_v:\Zb_v\to\Xb_v$ of
  $\gb$ to $\Zb_v$ is a (possibly disconnected) tamely ramified
  $m$-cyclic cover. If $v$ is an interior vertex, then $\gb_v$ is
  branched at most at the points $\xb_e\in\Xb_v$, with $s(e)=v$.  The
  canonical generator of inertia above $\xb_e$ (with respect to
  $\chib$) is $a_e$, where $a_e$ is as in Lemma \ref{lem2.1}.  In
  other words, $\gb_v:\Zb_v\to\Xb_v\cong\PP^1_k$ is the $m$-cyclic
  cover of type $(\xb_e;a_e)$.
  
  For $i\in I$, the cover $\Zb_i\to\Xb_i$ is ramified at the two
  points $\xb_i$ and $\xb_e$ (where $e$ is the unique edge with
  $s(e)=i$). The canonical generator of inertia above $\xb_i$ (resp.\ 
  $\xb_e$) is $a_i$ (resp.\ $a_{e}=m-a_i$). Finally, for $v\in B-I$,
  the cover $\Zb_v\to\Xb_v$ is totally disconnected, i.e.\ identifies
  $\Zb_v$ with $m$ disjoint copies of $\Xb_v$.
\end{prop}

\proof By the definition of admissibility, the covers
$\gb_v:\Zb_v\to\Xb_v$ are as in the statement of the proposition, for
{\em certain} integers $a_e$. Moreover, the integers $a_e$
verify Condition (i), (ii) and (iii) of Lemma \ref{lem2.1}, at least
modulo $m$. Now the statement of Lemma \ref{lem2.1} says that we may
choose the $a_e$'s such that they verify Condition (i), (ii) and (iii)
exactly. 
\Endproof

\subsection{The reduction of the $p$-cyclic cover $Y\to Z$}
\label{sec2.2}

We continue with the notation and assumptions of the previous
subsection. Recall from Section \ref{sec1.1} that the \'etale
$p$-cyclic subcover $Y\to Z$ of $f$ is given (birationally) by an
equation of the form $y^p=u$, where $u$ is a certain rational function
on $Z$, unique up to multiplication with a $p$th power. We may assume
that the $K$-model $Y_K\to Z_K$ (obtained from the stable model of
$f$) is defined by the same equation. By assumption, $Y_K\to Z_K$
extends to a finite map $Y_R\to Z_R$ between semistable $R$-curves,
where $Y_R$ is the stably pointed model of $Y_K$ and
$Z_R=Y_R/(\ZZ/p)$. Neglecting the $\ZZ/m$-action on $Z_R$, this is
precisely the situation studied in \cite{Saidi00} (see also
\cite{GreenMatignon99} and \cite{Yannick}).

Choose a vertex $v\in V$, and write $\Yb_v$ for the inverse
image of $\Zb_v$ in $\Yb$. We distinguish the following three cases
(compare with \cite{Saidi00}, \S 1.4 and \cite{Yannick}, \S 5.1):

\begin{itemize}
\item {\em (multiplicative reduction)} Suppose that the map
  $\Yb_v\to\Zb_v$ is inseparable (therefore, $v\in V'$ is an interior
  vertex, by Assumption \ref{weakass}). Suppose, moreover, that the
  restriction $\ub_v:=u|_{\Zb_v}$ of $u$ to $\Zb_v$ is not a $p$th
  power (in the function field of $\Zb_v$). Then the cover
  $\Yb_v\to\Zb_v$ is given (birationally) by the equation $y^p=\ub_v$
  and carries a natural structure of $\bmu_p$-torsor. Let
  $\omega_v:=d\ub_v/\ub_v$. The differential $\omega_v$ is not zero,
  regular on $\Zb_v$ and does not depend on the choice of the rational
  function $u$. One easily checks that
  $\phi_a^*\omega_v=\chib(a)\omega_v$, for $a\in\ZZ/m$. We write
  $\omega_v\in H^0(\Zb_v,\Omega^1)_{\chib}$. Furthermore, the zero
  locus of $\omega_v$ is contained in $ \gb^{-1}(\{\xb_e\mid
  s(e)=v\})$.
\item {\em (additive reduction)} Suppose that $\Yb_v\to\Zb_v$ is
  inseparable (hence $v\in V'$), and that the restriction of $u$ to
  $\Zb_v$ is a $p$th power. Then, in a neighborhood of any point
  $\zb\in\Zb_v$ on $Z_R$ and after multiplying $u$ with a suitable
  $p$th power, we can write $u=1 + \pi^pw_v$, such that $\pi\in R$,
  $0<v_R(\pi)<v_R(p)/(p-1)$ and such that $\wb_v:=w_v|_{\Zb_v}$ is not
  a $p$th power.  The restriction of $\Yb_v\to\Zb_v$ to the open
  subset $\Ub_v:=\Zb_v-\gb^{-1}(\{\xb_e\mid s(e)=v\})$ carries a
  natural structure of $\balpha_p$-torsor, locally given by the
  equation $\tilde{y}^p=\wb_v$. The differential $\omega_v:=d\wb_v$ is
  not zero, independent of all the choices we have made, and is
  regular on $\Ub_v$. Again, one easily checks that
  $\phi_a^*\omega_v=\chib(a)\omega_v$, for $a\in\ZZ/m$, i.e.\ 
  $\omega_v\in H^0(\Ub_v,\Omega^1)_{\chib}$. Furthermore, $\omega_v$
  has no zero in $\Ub_v$.
\item   
  {\em (\'etale reduction)} If $v$ is a leaf, then $\Yb_v\to\Zb_v$ is
  generically \'etale, ramified only at $\gb_v^{-1}(\xb_e)$, where $e$
  is the unique edge with $s(e)=v$.  Choose a point $\zb\in\Zb_v$
  above $\xb_e\in\Xb_v$.  In a neighborhood of $\zb$, the cover
  $\Yb_v\to\Zb_v$ is an Artin-Schreier cover, with equation
  $y^p-y=w^{-h_v}$, where $w$ is some local coordinate for $\Zb_v$ at
  $\zb$ and $h_v$ is the {\em conductor}.
\end{itemize}
We shall say that the vertex $v$ has multiplicative, additive or
\'etale reduction, according to which of the three cases occurs. 

The data $(\omega_v,h_v)$, which we obtain from the reduction of the
\'etale $p$-cyclic cover $Y\to Z$ satisfy certain compatibility
conditions, see \cite{Saidi00}. In our situation, they can be
formulated as follows. Let $e$ be an edge, and let
$\zb_e\in\Zb$ be a point above $\xb_e\in\Xb$. 
Define 
\begin{equation} \label{eq2.1a}
  h_e \; :=\; \begin{cases} 
     \;\;\ord_{\zb_e}(\omega_v)+1\   & \text{if $v:=s(e)\in V'$,}  \\
     \;\;\;\;-h_v                    & \text{if $v:=s(e)\in B$.} 
      \end{cases}
\end{equation}
Since $\omega$ is an eigenvector under the $\ZZ/m$-action, $h_e$ is
well defined. As a special case of \cite{Saidi00}, Cor.\ 2.8, we obtain 
\begin{equation} \label{eq2.1b}
    h_e + h_{\bar{e}} \;=\; 0,
\end{equation}
for all edges $e$. (It is also possible to derive \eqref{eq2.1b} from
\cite{Yannick}, Proposition 5.2.1, by considering the formal annulus
corresponding to the a point $\yb\in\Yb$ above $\xb_e\in\Xb$).

Using the fact that the data $(\omega_v;h_i)$ has to be compatible
with the $\ZZ/m$-action, we can express the numbers $h_e$ in terms of
certain (more convenient) numbers $\nu_e$. This generalizes
Proposition \ref{prop1.3}. 

\begin{prop} \label{prop2.2}
  For each edge $e$, there exists an
  integer $\nu_e$ such that 
  \begin{equation} \label{eq2.1c}
     h_e\;=\; \frac{\nu_em+a_e}{(a_e,m)}.
  \end{equation}
  The collection $(\nu_e)$ satisfies
  \begin{equation} \label{eq2.1d}
      \nu_e + \nu_{\bar{e}} \;=\; -1
  \end{equation}  
  and
  \begin{equation} \label{eq2.1e}
       \sum_{s(e)=v} (\nu_e-1) \;=\; -3,
  \end{equation}
  for all interior vertices $v$.
\end{prop}

\proof Recall that $a_e\in\ZZ/m$ is the canonical generator of inertia
for the $m$-cyclic cover $\gb_v:\Zb_v\to\Xb_v$ above $\xb_e$, see
Proposition \ref{prop2.1}. If $a_e\in\{0,m\}$ then $(a_e,m)=m$, and
the existence of an integer $\nu_e$ as in \eqref{eq2.1c} is trivial.
Now assume $0<a_e<m$, i.e.\ $\Zb_v\to\Xb_v$ is actually branched at
$\xb_e$, of order $m/(a_e,m)$. Then the vertex $s(e)$ is either an
interior vertex, or an element of $I$. Suppose that $v:=s(e)$ is an
interior vertex. Using $\phib_a^*\omega_v=\chib(a)\omega_v$ one shows
that $\ord_{\zb_e}(\omega_v)\equiv a_e/(a_e,m)-1\mod{m/(a_e,m)}$, see
\cite{IreneAux}, Lemma 2.6. This proves the existence of $\nu_e$ in
this case.  The case $i=s(e)\in I$ has already been proved in
Proposition \ref{prop1.3}. This completes the proof of the existence
of the integers $\nu_e$. Now \eqref{eq2.1d} follows from
\eqref{eq2.1b}, \eqref{eq2.1c} and the equality $a_e+a_{\bar{e}}=m$ by
a straightforward calculation. Finally, \eqref{eq2.1e} is a direct
consequence of Proposition \ref{prop2.1}, Equation \eqref{eq2.1c} and
the Riemann--Hurwitz formula.  \Endproof

We shall call an edge $e$ {\em terminal} if $v:=t(e)$ is a leaf; for
such an edge, $h_v=h_e=(\nu_em+a_e)/(a_e,m)$ is the conductor of the
Artin-Schreier cover $\Yb_v\to\Zb_v$. In particular, $h_v=\nu_e$ if
$v=t(e)\in B-I$. If $f$ is special, then $I=B$; moreover, for each
terminal edge $e$ (with $i:=t(e)\in I$), the integer $\nu_e=\nu_i$ is
either $0$ or $1$ and takes the value $0$ for exactly three terminal
edges (see Definition \ref{specialdef}). In this case, we can also say
a lot about the values $\nu_e$ on all edges $e$.

\begin{lem} \label{lem2.2}
  Assume that $f$ is special. Then:
  \begin{enumerate}
  \item 
    The integers $\nu_e$ (defined in Proposition \ref{prop2.2}) lie
    between $-2$ and $1$.
  \item
    There exists a unique interior vertex $v_0\in V'$ such that
    $\nu_e\geq 0$ for all edges $e$ with source $v_0$.
  \item
    If $v\neq v_0$ is an interior vertex, then there exists
    a unique edge $e$ with source $v$ such that $\nu_e<0$.
  \end{enumerate}
\end{lem}

\proof Since $f$ is special, we may assume that $\nu_i=0$ for
$i=1,2,3$ and $\nu_i=1$ for $i=4,\ldots,r$. Furthermore, $I=B$. For
any edge $e\in E$, let $I_e\subset I$ be the set of leaves $i\in I$
which `lie in the direction of $e$'. More precisely, $i\in I_e$ if and
only if $i$ lies in the connected component of $T-\{e\}$ which
contains the vertex $t(e)$. We claim that
\begin{equation} \label{eq2.2}
  \nu_e \;=\; 1 \;-\; |\,I_e\cap\{1,2,3\}\,|,
\end{equation}
for all $e\in E$. Let us check that the lemma follows from this claim.
Indeed, (i) is a trivial consequence of \eqref{eq2.2}, and part (ii)
and (iii) of the lemma follow once we observe that $v_0\in V$ has to
be the {\em median} of the three leaves $i=1,2,3$.

Recall that $\nu_e=\nu_i$ if $i:=t(e)\in I$. Using $\nu_i=0$ for
$i=1,2,3$ and $\nu_i=1$ for $i>3$, we conclude that \eqref{eq2.2}
holds for all edges $e$ such that $t(e)\in I$. For a general edge
$e\in E$, define $\nu_e':=1 \;-\; |\,I_e\cap\{1,2,3\}\,|$.  An easy
verification shows that the function $e\mapsto \nu_e'$ verifies
Equations \eqref{eq2.1d} and \eqref{eq2.1e}. We conclude that
$\nu_e=\nu_e'$ for all edges $e\in E$, by induction. This finishes the
proof of the lemma.  \Endproof

\begin{rem} \label{rem2.1}
  It is shown in \cite{Yannick} that the integers $h_e$ determine the
  radii of the formal annuli corresponding to the singular points of
  $Y_R$. To be more precise, let $\yb$ be an ordinary double point of
  the special fiber of a semistable $R$-curve $Y_R$. Then the complete
  local ring of $Y_R$ at $\yb$ is of the form $\Od_{Y_R,\yb}\cong
  R[[u,v\mid uv=\pi]]$, with $\pi\in R$. We define the {\em thickness}
  of $Y_R$ at $\yb$ as the (positive rational) number
  $\epsilon(Y_R,\yb):=v_R(\pi)$ (recall that $v_R(p)=1$).  Suppose
  $\yb\in\Yb$ is a point above $\xb_e\in\Xb$, $e\in E$. Suppose,
  moreover, that $v:=s(e)$ is a vertex with multiplicative reduction
  and $v':=t(e)$ has additive reduction. Then
  \begin{equation} \label{thicknesseq1}
       0 \;<\; \epsilon(Y_R,\yb) \;<\; \frac{1}{(p-1)\,h_e}.
  \end{equation}
  On the other hand, if $v:=s(e)$ has multiplicative and $v':=t(e)$
  \'etale reduction, then
  \begin{equation} \label{thicknesseq2}
       \epsilon(Y_R,\yb) \;=\; \frac{1}{(p-1)\,h_e}.
  \end{equation}
  This follows immediately from \cite{Yannick}, Chap. 5, 
  Proposition 2.1. Moreover, using \cite{Raynaud98}, Proposition
  2.3.2, one shows that
  \begin{equation} \label{thicknesseq3}
      \epsilon(X_R,\xb_e) \;=\; 
        \frac{p\,a_e}{(a_e,m)}\cdot\epsilon(Y_R,\yb).
  \end{equation}
\end{rem}

\subsection{The proof of Theorem \ref{thm1}}  \label{sec2.3}

Let $v_0$ be the `median vertex' of Lemma \ref{lem2.2} (ii).  Theorem
\ref{thm1} is equivalent to the statement that $v_0$ is the unique
interior vertex. Therefore, let us assume that there exists another
interior vertex $v\neq v_0$, and then try to arrive at a
contradiction.  By Lemma \ref{lem2.2} (iii), there exists a unique
edge $e$ with $s(e)=v$ such that $\nu_e<0$. This means that the
differential $\omega_v$ has a pole in each point $\zb_e\in\Zb_v$ above
$\xb_e\in\Xb_v$. If the cover $Y\to Z$ had multiplicative reduction at
the component $\Zb_v$, then the differential $\omega_v$ would be
regular on $\Zb_v$. Therefore, we have additive reduction at $\Zb_v$.
In particular, the differential $\omega_v$ is a nonzero {\em exact}
differential, i.e.\ of the form $\omega_v=du$, for some rational
function $u$ on $\Zb_v$.  Moreover, the divisor $(\omega_v)$ is
completely determined by the integers $\nu_e$, where $e$ runs through
the set of edges with source $v$.  By Lemma \ref{lem2.2} (i) and
Proposition \ref{prop2.2}, these numbers satisfy $-2\leq\nu_e\leq 1$
and $\sum_{s(e)=v}(\nu_e-1)=-3$. The following lemma gives the desired
contradiction, and thus finishes the proof of Theorem \ref{thm1}.

\begin{lem}  \label{fundlem}
  Let $k$ be an algebraically closed field of characteristic $p>0$,
  $m>1$ an integer dividing $p-1$, and $r\geq 3$.  Let $g:Z\to\PP^1_k$
  be an $m$-cyclic cover given (birationally) by an equation of the form
  \[
        z^m \;=\; \prod_{i=1}^r \;(x-x_i)^{a_i},
  \]
  with $x_1,\ldots,x_r\in k$ pairwise distinct, and $0<a_i<m$ such
  that $\sum_i a_i=m$. Let $\phi:Z\iso Z$ be a generator of
  $\Aut(Z/\PP^1_k)$, such that $\phi^*z=\zeta z$, $\zeta\in k^\times$
  an $m$th root of unity. For $i=1,\ldots,m$, let $m_i:=m/(a_i,m)$,
  $\at_i:=a_i/(a_i,m)$ and $P_i:=g^{-1}(x_i)$ (considered as a divisor
  on $Z$). Let $\omega$ be a meromorphic differential form on $Z$ such
  that
  \begin{enumerate}
  \item
    $\phi^*\omega=\zeta\omega$, and 
  \item 
    $(\omega) \;=\; \sum_{i=1}^r\; (m_i\nu_i+\at_i)P_i$, with
    integers $-2\leq\nu_i\leq 1$ such that $\sum_i\nu_i=r-3$,
    $\nu_1<0$ and $\nu_i\geq 0$ for $i>1$.
  \end{enumerate}
  Then $\C(\omega)\neq 0$, i.e.\ $\omega$ is not exact. 
\end{lem}
 
\proof Assume that $\omega$ is a meromorphic differential on $Z$
such that (i) and (ii) hold. After a change of coordinate, we may
assume that $x_1=0$. By (i), we can write $\omega=f\,z\,{\rm d}x$,
where $f$ is a rational function in $x$. Expanding $f$ as a
Taylor series at $x=x_1=0$, we obtain
\begin{equation} \label{eq2.3}
   \omega \;=\;  (\sum_{j=-3}^\infty\, c_j\,x^j)\,z\,{\rm d}x,
\end{equation}
with $c_j\in k$. Note that $x$ (resp.\ $z$) has a zero of order $m_1$
(resp.\ of order $\at_1$) at each point $\zb\in P_1$. In particular,
the coefficients $c_{-3}$ and $c_{-2}$ contribute to the poles of
$\omega$ in $P_1$, which are of order $(-m_1\nu_1+a_1-1)$, by (ii). 

{\bf Claim 1:} There exist elements $b_1,b_2\in k$ such that
\begin{equation} \label{eq2.4}
  \omega':=\omega-{\rm d}u\;\;\; \text{\rm is regular on $Z$, where}\;\;\;\;
  u:=(b_1x^{-2}+b_2x^{-1})\,z.
\end{equation}

Since $\omega$ and $u$ are regular away from $P_1$, we only
have to pay attention to the points in $P_1$.
We compute `Taylor series' as in \eqref{eq2.3}:
\begin{equation} \label{eq2.5}
   {\rm d}z \;=\; (\frac{a_1}{m}\,x^{-1} + d_0 + d_1\,x
   +\cdots)\,z\,{\rm d}x
\end{equation}
and
\begin{equation} \label{eq2.6}
\begin{split}
   {\rm d}u \;&=\; (-2b_1\,x^{-3}-b_2\,x^{-2})\,z\,{\rm d}x \;+\;
             (b_1\,x^{-2}+b_2\,x^{-1})\,{\rm d}z            \\
        &=\; \Big( (\frac{a_1}{m}-2)\,b_1\,x^{-3} \;+\;
                   (d_0\,b_1+(\frac{a_1}{m}-1)\,b_2)\,x^{-2}
                   \;+\; \cdots \Big)\,z\,{\rm d}x.              \\
\end{split}
\end{equation}                      
Hence, to prove Claim 1, we have to find $b_1,b_2$ such that
\begin{equation} \label{eq2.7}
  (\frac{a_1}{m}-2)\,b_1 \;=\; c_{-3}, \qquad 
  d_0\,b_1  +\; (\frac{a_1}{m}-1)\,b_2  \;=\; c_{-2}.
\end{equation}
Using $m|p-1$ and $a_1\leq m-2$ one shows that $p$ does not divide
$2m-a_1$ and $m-a_1$; therefore, $a_1/m-2,a_1/m-1\neq 0$ in $k$, and
we can solve \eqref{eq2.7} in $b_1$ and $b_2$. This proves Claim 1.

{\bf Claim 2:} $\omega \neq {\rm d}u$. \\
Assuming the contrary, we would have
\begin{equation} \label{eq2.8}
   (\omega) \;=\; ({\rm d}u) \;=\; \sum_{i=1}^r \; (m_i\nu_i+\at_i-1)P_i,
\end{equation}
by Condition (ii). Assume for the moment that the order of $u$ at all
the ramification points is prime-to-$p$. Then \eqref{eq2.8} implies
\begin{equation} \label{eq2.9}
   (u) \;\geq\; \sum_i\; (m_i\nu_i+\at_i)P_i.
\end{equation}
But the divisor on the right has degree $\sum_i m\nu_i+a_i=m(r-2)>0$,
contradiction! Thus, in order to prove Claim 2, it suffices to show
that $u$ has no zero or pole in one of the ramification points of
order divisible by $p$. Let $z_i\in P_i$ be a ramification point above
$x_i$, for some $i$. Since $\phi^*u=\zeta u$,
$\ord_{z_i}(u)=m_ik+\at_i$, for some integer $k$. If $p$ divides
$\ord_{z_i}(u)=m_ik+\at_i$ then either $k\geq 1$ or $k<0$. In the
first case, we would have $\ord_{z_i}({\rm
d}u)=m_i\nu_i+\at_i>m_i+\at_i$, hence $\nu_i>1$, which contradicts our
assumptions. The second case can occur only for $i=1$.  But for $i=1$,
$k\in\{-2,-1\}$, and in this case we have already shown that
$km_i+\at_i$ is prime to $p$ (see the end of the proof of Claim 1). We
conclude that $\omega\neq{\rm d}u$, as asserted by Claim 1.

Set $\omega':=\omega-{\rm d}u$, and note that
$\phi^*\omega'=\zeta\omega'$, by Condition (i) and the definition of
$u$. Following our previous notation, we can write $\omega'\in
H^0(Z,\Omega^1_{Z/k})_\chi$, where $\chi:\ZZ/m\to k^\times$ is the
character with $\chi(1)=\zeta$. By Lemma \ref{lem1.0} (ii), we have
$\C(\omega)=\C(\omega')\neq 0$. This proves the lemma, and also
Theorem \ref{thm1}.  \Endproof

\subsection{The monodromy group of a special cover} \label{sec2.4}

The analysis of the stable reduction of a special cover shows
somewhat more than what is stated in Theorem \ref{thm1}. 
We use the same notation as in Theorem \ref{thm1}. In particular,
$f:Y\to\PP^1$ is a special cover of type $(x_i;a_i)$, with stable
reduction $\fb:\Yb\to\Xb$. Define
\begin{equation}
  \D_i \;:=\; \{\; x\in\PP^1(\Kb) \mid\;
     \text{$x$ specializes to a point on 
                 $\Xb_i-\Xb_0\cong\AA^1_k$}\;\},
\end{equation}
the closed rigid disk containing all points of $\PP^1$ which
specialize to $\Xb_i$. In particular, $x_i\in\D_i$.  

\begin{prop}  \label{prop2.3}
  We have 
  \[
     \D_i \;=\; \{\; x\in \PP^1(\Kb) \;\mid\; 
          v_R(x-x_i) \geq \frac{p\,m_i}{(p-1)\,h_i} \;\}.
  \]
  (Recall that $m_i=m/(a_i,m)$ and $h_i=(m\nu_i+a_i)/(a_i,m)$.)
\end{prop}

\proof
This follows immediately from Remark \ref{rem2.1}, Equations
\eqref{thicknesseq2} and \eqref{thicknesseq3}. 
\Endproof
   
For the rest of this section, we assume that the absolute ramification
index of $K_0$ is one. Thus, we can identify $R_0$ with the ring
$W(k)$ of Witt vectors over the residue field $k$ (in view of the
results of Section \ref{sec3}, this is not a serious restriction). By
Proposition \ref{prop1.2}, the field of moduli of $f$ is
$K\inn=K_0(\zeta_p^{(m)})$. Let $K=K\stab$ be the minimal extension of
$K\inn$ over which $f$ has stable reduction. One can show the
following, see \cite{Raynaud98} and \cite{badred}. The extension
$K/K\inn$ is Galois, of degree prime-to-$p$. The Galois group
$\Gamma:=\Gal(K/K\inn)$ acts faithfully and $k$-linearly on $\Yb$
(where $\fb:\Yb\to\Xb$ is the stable reduction of $f$), and this
action commutes with the action of $H$. Therefore, we get an induced
action of $\Gamma$ on $\Xb$. The group $\Gamma$ is called the {\em
  monodromy group} of $f$.

\begin{thm} \label{thm1b}
  Let $f:Y\to\PP^1$ be a special cover of type $(x_i;a_i)$. Assume
  that the branch points $x_i$ are rational over $K_0$, the fraction
  field of $W(k)$. Then the order of the monodromy group $\Gamma$ of
  $f$ is
  \[
      |\Gamma| \;=\; [K:K\inn] \;=\; m\cdot{\rm lcm}(h_1,\ldots,h_r).
  \]
  Furthermore, the action of $\Gamma$ is trivial on $\Xb_0$ and cyclic
  of order $h_i(p-1)/m_i$ on $\Xb_i$.
\end{thm}

\proof The proof of Theorem \ref{thm1} shows that here exists an open
subset $U\subset Z_R$ such that $U\otimes_R k\subset\Zb_0$ is nonempty
and $V:=U\times_{Z_R}Y_R\to U$ is a $\bmu_p$-torsor. But the generic
fiber $V_K\to U_K$ is a $\ZZ/p$-cover. Therefore, the extension
$K/K_0$ contains the $p$th roots of unity. Moreover, the subgroup
$\Gammat:=\Gal(K/K_0(\zeta_p))\subset\Gamma$ is precisely the
stabilizer of $\Yb_0\subset\Yb$ (where $\Yb_0$ denotes the inverse
image of $\Zb_0$). It follows from Proposition \ref{prop1.2} (ii) that
$\Gamma/\Gammat\cong\ZZ/m$. Recall that $\Yb$ is the union of $\Yb_0$
and $\Yb_i$, the (possibly disconnected) inverse image of $\Xb_i$, for
$i=1,\ldots,r$. Let $\yb$ be a point where $\Yb_i$ intersects $\Yb_0$,
and let $\Yb_i'$ be the connected component of $\Yb_i$ containing
$\yb$. Since $\Gammat$ acts trivially on $\Yb_0$, it fixes $\Yb_i'$
and $\yb$. By Remark \ref{rem2.1}, Equation \ref{thicknesseq3},
$\epsilon(Y_R,\yb)=1/(p-1)h_i$. This means that the complete local
ring of $Y_R$ at $\yb$ is of the form $R[[u,v\mid
uv=\lambda^{1/h_i}]]$, with $\lambda:=\zeta_p-1\in K_0(\zeta_p)$.
Applying an element of $\Gammat$ to the equation $uv=\lambda^{1/h_i}$,
one shows that $\Gammat$ induces a cyclic action on $\Yb_i'$, of order
$h_i$. More precisely, the image of $\Gammat$ in $\Aut_k(\Yb_i')$ is
the quotient
$\Gal(K_0(\zeta_p,\lambda^{1/h_i})/K_0(\zeta_p))\cong\ZZ/h_i$. Note
that the action of $\Gammat$ on $\Yb_i'$ commutes with the action of
the decomposition group $H_i\subset H$ of $\Yb_i'$, which is of order
$pm_i$. Now the statement of Theorem \ref{thm1b} on the order of
$\Gamma$ follows from the fact that the action of $\Gamma$ on $\Yb$ is
faithful. The statement about the action of $\Gamma$ on $\Xb_i$
follows as well, using the fact that $h_i$ is relatively prime to
$pm_i$ (it can also be deduced directly from Proposition
\ref{prop2.3}).  \Endproof

%-------------------------------------------------------------------

\section{Construction of special covers} \label{sec3}

This section is concerned with the construction of special covers by
lifting certain objects from characteristic $p$ to characteristic $0$.
We start by defining {\em special degeneration data}, which are
essentially given by an $m$-cyclic cover $\Zb_0\to\PP^1_k$ of the
projective line in characteristic $p$, together with a logarithmic
differential form $\omega_0$ on $\Zb_0$, with certain prescribed
zeros. It is immediate from the results of the previous section that
the reduction $\fb:\Yb\to\Xb$ of a special cover $f$ corresponds
essentially to a special degeneration datum. The main result of this
section (Theorem \ref{thm2}) states that, conversely, every special
degeneration datum arises as the reduction of a special cover $f$.
Moreover, the cover $f$ is essentially unique, once we have chosen the
branch points $x_i$. 

The proof of Theorem \ref{thm2} is divided into two steps. In the
first step, we lift the $\bmu_p$-torsor $\Yb_0\to\Zb_0$ corresponding
to the differential $\omega_0$ to characteristic $0$, in a
$\ZZ/m$-equivariant way. This construction yields a metacyclic cover
$f:Y\to\PP^1$, which is essentially unique because the
$\chib$-eigenspace of the $p$-torsion of the Jacobian of $\Zb_0$ is
\'etale (Proposition \ref{prop1.2}). In the second step we show
that the cover $f$ we have constructed is special provided that we
have chosen the branch points $x_i$ inside certain closed rigid disks
$\D_i\subset\PP^1$. The proof uses the {\em monodromy action} on the
stable reduction of Galois covers, and a deformation argument.

In Section \ref{sec3.5} we determine all special degeneration data in
the case $r=4$.

\subsection{Special degeneration data} \label{sec3.1}

Let $f:Y\to\PP^1$ be a special cover of type $(x_i;a_i)$, with stable
model $f_R:Y_R\to X_R$ and reduction $\fb:\Yb\to\Xb$ (see Definition
\ref{specialdef}). By Theorem \ref{thm1}, $\Xb$ consists of $r+1$
components $\Xb_0,\ldots,\Xb_r$, such that, for $i\geq 1$, $\Xb_i$ is
the tail containing the specialization $\xb_i$ of the branch point
$x_i$. The component $\Xb_0$, which intersects all of the
components $\Xb_i$, $i\geq 1$, is called the {\em original component}
of $\Xb$.  We have a canonical isomorphism $\Xb_0\cong\PP^1_k$ arising
from the contraction morphism $q:X_R\to\PP^1_R$. This isomorphism
identifies the intersection point $\tau_i\in\Xb_0\cap\Xb_i$ with the
specialization of $x_i$, regarded as point on $\PP^1_R$.  We may 
assume, without loss of generality, that $\tau_i\neq\infty$.

We have seen in Section \ref{sec2.1} that the map $\fb:\Yb\to\Xb$ is
the composition of an $m$-cyclic admissible cover $\gb:\Zb\to\Xb$ with
a finite map $\Yb\to\Zb$ of degree $p$ which is the reduction of the
\'etale $p$-cyclic cover $Y\to Z$.  By Proposition \ref{prop2.1}, the
restriction $\gb_0:\Zb_0\to\Xb_0=\PP^1_k$ of $\gb$ to the original
component can be identified with the $m$-cyclic cover of type
$(\tau_i;a_i)$. Moreover, the $p$-cyclic cover $Y\to Z$ has
multiplicative reduction at $\Zb_0$.  This means that the induced
cover $\Yb_0\to\Zb_0$ carries the structure of a $\bmu_p$-torsor. This
structure gives rise to a regular differential form $\omega_0$ such
that $\phib_a^*\omega_0=\chib(a)\omega_0$ (recall that
$\phib_a:\Zb\iso\Zb$ is the automorphism induced by $a\in\ZZ/m$; we
write $\omega_0\in H^0(\Zb_0,\Omega_1)_{\chib}$). Let
$m_i:=m/(a_i,m)$, $\at_i:=a_i/(a_i,m)$ and $P_i:=\gb_0^{-1}(\tau_i)$
(we regard $P_i$ as a divisor on $\Zb_0$). By Proposition
\ref{prop2.2}, there exist integers $\nu_i\in\{0,1\}$ with
$\sum_i\nu_i=r-3$ such that
\begin{equation} \label{eq3.1}
    (\omega_0) \;=\; \sum_i\,(m_i\nu_i+\at_i-1)\,P_i.
\end{equation}
Furthermore, $\omega_0$ is logarithmic. In terms of the Cartier
operator $\C$, this means that 
\begin{equation} \label{eq3.2}
   \C(\omega_0) \;=\; \omega_0.
\end{equation}

\begin{defn} \label{degendef}
  Let $k$ be an algebraically closed field of characteristic $p>0$. A
  {\em special degeneration datum} over $k$ is given by
  \begin{itemize}
  \item 
    pairwise distinct $k$-rational
    points $\tau_1,\ldots,\tau_r\in\PP^1_k$, with $r\geq 3$,
  \item 
    an integer $m>1$ dividing $p-1$, and integers $a_1,\ldots,a_r$
    such that $0<a_i<m$ and $\sum_i\,a_i=m$, (we let
    $\gb_0:\Zb_0\to\PP^1_k$ be the $m$-cyclic cover of type
    $(\tau_i;a_i)$; furthermore, we set $m_i:=m/(a_i,m)$,
    $\at_i:=a_i/(a_i,m)$ and $P_i:=\gb_0^{-1}(\tau_i)$),
  \item
    integers $\nu_1,\ldots,\nu_r\in\{0,1\}$ such that 
    $\sum_i\nu_i=r-3$, and 
  \item
    a differential form $\omega_0\in H^0(\Zb_0,\Omega^1)_{\chib}$,
    such that \eqref{eq3.1} and \eqref{eq3.2} hold.
  \end{itemize}
\end{defn}

As explained in the paragraph preceding Definition \ref{degendef}, we
can attach to (the reduction of) a special cover $f:Y\to\PP^1$ of type
$(x_i;a_i)$ a special degeneration datum $(\tau_i;a_i;\nu_i;\omega_0)$.
Theorem \ref{thm2} below states that, conversely, every special
degeneration datum arises in this way.

For the rest of this section, we fix a special degeneration datum
$(\tau_i;a_i;\nu_i;\omega_0)$ over $k$. Let $K_0$ denote the fraction
field of $R_0:=W(k)$, the ring of Witt vectors over $k$. Choose an
algebraic closure $\Kb$ of $K_0$. For $i\in I$, choose a
$K_0$-rational point $\xt_i\in\PP^1(K_0)$ which lifts
$\tau_i\in\PP^1(k)$. Let
\begin{equation} \label{diskeq}
   \D_i \;:=\; \{\; x\in \PP^1(\Kb) \;\mid\; v_R(x-\xt_i) 
         \geq \frac{p\,m_i}{(p-1)\,h_i} \;\}
\end{equation}
(compare to the statement of Theorem \ref{thm1b}). We claim that the
collection of disks $(\D_i)$ does not depend, up to an automorphism of
$\PP^1_{K_0}$, on the choice of the points $\xt_i$. To show that this
is so, we may assume that $\nu_1=\nu_2=\nu_3=0$. Furthermore, we can
always normalize our choice such that $\xt_1=0$, $\xt_2=1$ and
$\xt_3=\infty$. For $i=4,\ldots,r$, we have $0<pm_i/(p-1)h_i<1$. Using
the triangle inequality and the fact that the valuation $v_R$ takes
integral values on $K_0$, one shows that $\D_i$ does not depend on the
choice of $\xt_i$, for $i=4,\ldots,r$.

\begin{thm}  \label{thm2}
  Let $(\tau_i;a_i;\nu_i;\omega_0)$ be a special degeneration datum over
  $k$ and let $x_1,\ldots,x_r$ be $\Kb$-rational points on $\PP^1$,
  such that $x_i\in\D_i$. Then there exists a special cover
  $f:Y\to\PP^1$ of type $(x_i;a_i)$, unique up to isomorphism, which
  gives rise to $(\tau_i;a_i;\nu_i;\omega_0)$.
\end{thm}

This theorem corresponds to Theorem B in the introduction. The proof
is given in the next three subsections. It seems reasonable to expect
that the condition $x_i\in\D_i$ is also necessary for $f$ to be
special. Unfortunately, the method of proof we use here does not give
such an `if and only if' result.

\subsection{Definition of the lift} \label{sec3.2}

Let $(\tau_i;a_i;\nu_i;\omega_0)$ be a special degeneration datum over
$k$ and $\gb_0:\Zb_0\to\PP^1_k$ the $m$-cyclic cover of type
$(\tau_i;a_i)$. Choose $\Kb$-rational points $x_i\in\PP^1(\Kb)$
lifting the points $\tau_i$ (for the moment, we do not assume that
$x_i\in\D_i$). Choose a finite extension $K/K_0$ such that $x_i$ is
$K$-rational; let $R$ be the ring of integers of $K$.

The $m$-cyclic cover $\gb_0$ lifts uniquely to an $m$-cyclic cover
$g_{R}':Z_{R}'\to\PP^1_{R}$ of smooth curves, tamely ramified along
the closure of $\{x_1,\ldots,x_r\}\subset\PP^1_{K}$ inside
$\PP^1_{R}$. Let $g_{K}:Z_{K}\to\PP^1_{K}$ be the generic fiber of
$g_{R}'$ and $g:Z\to\PP^1$ its base change to $\Kb$. Note that $g$ is
the $m$-cyclic cover of type $(x_i;a_i)$. Let $J_{R}$ be the
N\'eron-model of $J_{Z_{K}}$ over $R$. Since $Z_{R}'$ is smooth over
$R$, $J_{R}$ represents the functor $\Pic^0(Z_{R}'/R)$, see
\cite{BoLuRa}, \S 9, Proposition 4. The universal property of the
N\'eron model defines a surjective specialization map $J_{Z_{K}}(K)\to
J_{\Zb_0}(k)$. As we have seen in Section \ref{sec1.2b}, the
specialization map induces an isomorphism
\begin{equation} \label{eq3.4}
  J_Z[p]_{\chib} \liso J_{\Zb_0}[p]_{\chib}
\end{equation}
of $\FF_p$-modules (of rank $r-2$). 

The logarithmic differential $\omega_0$ corresponds to a line bundle
$\Lb$ on $\Zb_0$, in the following way (see e.g.\ \cite{Milne}, III,
\S 4). Let $\ub$ be a rational function on $\Zb_0$ such that
$\omega_0=d\ub/\ub$. Then $(\ub)=p\cdot \Db$, for a divisor $\Db$ of
degree $0$ on $\Zb_0$; we set $\Lb:=\O_{\Zb_0}(\Db)$. By definition,
$\Lb^{\otimes p}\cong\O_{\Zb_0}$.  Moreover,
$\phib_a^*\omega_0=\chib(a)\omega_0$ implies
$\phib_a^*\Lb\cong\Lb^{\otimes \chib(a)}$ (note that this makes sense
because $\chib(a)\in\FF_p^\times$). In other words, $\Lb$ corresponds
to an element of $J_{\Zb_0}[p]_\chi$. Let $L$ be the line bundle on
$Z$ corresponding to $\Lb$ under the isomorphism \eqref{eq3.4}. By the
definition of the specialization map \eqref{eq3.4}, $L$ is actually
the pullback of a line bundle $L_{R}$ on $Z_{R}'$, and $L_{R}\cong
\O_{Z_{R}'}(D)$, where $D$ is a horizontal divisor on $Z_{R}'$ such
that $p\cdot D=(u)$ for some rational function $u$. By construction,
we have $\omega_0=d\ub/\ub$, with $\ub:=u|_{\Zb_0}$.  We let $Y\to Z$
be the $\bmu_p$-torsor corresponding to $L$ (birationally given by the
equation $y^p=u$).  After choosing a $p$th root of unity
$\zeta_p\in\Kb$, we can regard $Y\to Z$ as an \'etale $p$-cyclic
cover. Now the composition $f:Y\to Z\to\PP^1$ is a metacyclic cover of
type $(x_i;a_i)$ (see Section \ref{sec1.1}).

The cover $f$ we have constructed will not be special, in general.
However, if $f$ is special then, by construction, it gives rise to the
special degeneration datum $(\tau_i;a_i;\nu_i;\omega_0)$ we started
with. It is also clear that any special cover which gives rise to
$(\tau_i;a_i;\nu_i;\omega_0)$ is isomorphic to $f$. Thus, in order to
prove Theorem \ref{thm2}, we have to show that $f$ is special provided
that $x_i\in\D_i$, for all $i\in I$. Before we can give a proof of
this claim (in Section \ref{sec3.4}), we need to analyze the stable
reduction of $f$. For this step, it is not yet necessary to assume
$x_i\in\D_i$.

\subsection{Analyzing the stable reduction of $f$} \label{sec3.3}

We may assume that the cover $f$ constructed above has stable
reduction over the field $K$. Let $f_R:Y_R\to X_R$ be the stable model
of $f$, and $\fb:\Yb\to\Xb$ its reduction.  For $i=1,\ldots,r$, let
$\xb_i\in\Xb$ be the specialization of the branch point $x_i$, and let
$\Xb_i$ be the component of $\Xb$ containing $\xb_i$. Since
$x_i\not\equiv x_j\mod{v_{R_0}}$ for $i\neq j$ (as points on $\PP^1$),
it follows from \cite{Raynaud98}, \S 3, that the components $\Xb_i$
are pairwise distinct tails of $\Xb$ and that $\fb:\Yb\to\Xb$ is
separable exactly over the tails. In other words, Assumption
\ref{weakass} holds.  Note that the stronger Assumption
\ref{specialass} may not hold, as $\Xb$ might have new tails.

However, Assumption \ref{weakass} being valid, we may use the notation
set up in Section \ref{sec2.1} and \ref{sec2.2}, concerning the
structure of $\fb:\Yb\to\Xb$ as the composition of the $m$-cyclic
admissible cover $\gb:\Zb\to\Xb$ and the `mixed torsor' $\Yb\to\Zb$.
Recall that we described this structure using certain combinatorial
data $(T;a_e;\nu_e)$. Here $T$ is the dual graph of components of the
semistable curve $\Xb$. The integers $a_e$ (where $e$ is an edge of
$T$) describe the admissible $m$-cyclic cover $\gb:\Zb\to\Xb$.
Finally, the integers $\nu_e$ (together with the $a_e$) determine the
order of the zeros and the poles of the differentials $\omega_v$
attached to the torsors $\Yb_v\to\Zb_v$, where $v$ is an interior
vertex of $T$.  By construction of $f$, $Y\to Z$ has multiplicative
reduction above the original component $\Xb_0$, and the resulting
$\bmu_p$-torsor $\Yb_0\to\Zb_0$ corresponds to the differential
$\omega_0$. Therefore, \eqref{eq3.1} implies that
\begin{equation} \label{eq3.4a}
  \nu_{e_i} \;=\; \nu_i\in\{0,1\},
\end{equation}
where $\nu_i$ is given by the special degeneration datum and $e_i$ is
the edge corresponding to the point $\tau_i\in\Xb_0$ (in particular,
$\tau_i$ is a singular point of $\Xb$). Note that we do not know
(unless $f$ is special) whether $h_i=\nu_im_i+\at_i$ is the conductor 
of the Artin-Schreier cover $\Yb_i\to\Zb_i$ over the tail $\Xb_i$. 

Let $v_0\in V$ be the vertex of the tree $T$ corresponding to the
original component $\Xb_0$. For any edge $e$ of $T$, we let $T_e$ be
the connected component of $T-\{e\}$ which contains $t(e)$. We shall
call an edge $e$ {\em positive} if $v_0\not\in T_e$, i.e.\ if $e$ is
directed away from the vertex $v_0$.

\begin{lem} \label{lem3.1}
  \begin{enumerate} 
  \item 
    Let $e$ be a positive edge. Then $\nu_e\geq 0$. If, moreover,
    $a_e\equiv 0\mod{m}$, then $\nu_e>1$.
  \item 
    The points $\tau_1,\ldots,\tau_r$ are precisely the points of
    $\Xb_0$ which are singular points of $\Xb$.
  \item 
    Fix $i\in I=\{1,\ldots,r\}$, and let $e_i$ be the edge with
    source $v_0$ corresponding to $\tau_i$. If $\nu_i=0$, then
    $T_{e_i}=\{i\}$.  On the other hand, if $\nu_i=1$, then either
    $T_{e_i}=\{i\}$ or we are in the following case. The vertex
    $v:=t(e_i)$ is the source of exactly three edges,
    $\bar{e}_i,e',e''$. Also, $t(e')=i\in I$ and $t(e'')\in B-I$ is a
    leaf. See Figure \ref{treespic}.
  \end{enumerate}    
\end{lem}

\begin{figure}[bt]
\begin{center}
    \setlength{\unitlength}{0.0005in}
{\renewcommand{\dashlinestretch}{30}
\begin{picture}(7563,3321)(0,-10)
\put(801,462){\circle*{90}}
\put(171,2712){\circle{90}}
\put(1431,2712){\circle{90}}
\put(801,1722){\circle*{90}}
\put(5301,2532){\circle*{90}}
\drawline(801,1722)(216,2667)
\drawline(304.671,2580.759)(216.000,2667.000)(253.655,2549.178)
\drawline(801,1722)(1386,2667)
\drawline(1348.345,2549.178)(1386.000,2667.000)(1297.329,2580.759)
\drawline(801,462)(801,1677)
\drawline(831.000,1557.000)(801.000,1677.000)(771.000,1557.000)
\drawline(3951,912)(7551,912)           %\Xb_0
\drawline(6651,462)(6651,2892)          %\Xb_v
\drawline(7101,1812)(4851,1812)         %\Xb_i
\drawline(7101,2532)(4851,2532)         %\Xb_b
\dottedline{45}(801,462)(351,12)
\dottedline{45}(801,462)(1251,12)
\dottedline{45}(3951,912)(3051,912)
\put(351,462){$v_0$}
\put(1071,1092){$e_i$}
\put(150,2037){$e''$}
\put(1250,2037){$e'$}
\put(1431,3027){$i$}
\put(4370,2450){$\Xb_i$}
\put(6550,3162){$\Xb_v$}
\put(4176,500){$\Xb_0$}
\put(5250,2757){$\xb_i$}
\put(8500,1400){$\Xb$}
\put(-500,1400){$T$}
\put(6850,1900){$\tau'$}
\put(6850,2600){$\tau''$}
\put(6850,1000){$\tau_i$}
\end{picture}
}

%%% Local Variables: 
%%% mode: latex
%%% TeX-master: "special"
%%% End: 
\end{center}
\caption{}
\label{treespic}
\end{figure}

\proof Suppose we have a positive edge $e$ such that $a_e\equiv
0\mod{m}$ and $\nu_e\leq 1$.  Note that this implies $a_{e'}\equiv
0\mod{m}$ for all edges contained in the subtree $T_e$ (otherwise,
$T_e$ would contain exactly one leaf $i\in I$, and then $a_e\equiv
a_i\not\equiv 0\mod{m}$). Assume first that $v:=t(e)$ is not a leaf.
From \eqref{eq2.1d} and \eqref{eq2.1e} we deduce the inequality
\begin{equation} \label{eq3.5}
  \sum_{s(e')=v,\,e'\not=\bar{e}}(\nu_{e'}-1) \;=\;-1+\nu_e \;\leq\; 0.
\end{equation}
Thus, we have $\nu_{e'}\leq 1$ and $a_{e'}\equiv 0\mod{m}$ for at
least one positive edge $e'$ with source $v$. Hence, after a finite
number of steps, we find an edge $e$ such that $v:=t(e)$ is a leaf,
$a_e\equiv 0\mod{m}$ and $\nu_e\leq 1$. This means that $\Xb_v$ is a
new tail of $\Xb$ and that the conductor of the Artin-Schreier cover
$\Yb_v\to\Zb_v$ is $h_v=\nu_e\leq 1$. It follows that each connected
component of $\Yb_v$ is a tail of $\Yb$ of genus $0$; furthermore, no
ramification point specializes to $\Yb_v$. This contradicts the
minimality of the stable model $Y_R$, and proves the second assertion
of (i). The proof of the first assertion is similar, and uses the fact
that the conductor of an Artin-Schreier cover is $\geq 1$.

Statement (ii) of the lemma follows immediately from (i). Indeed, a
singular point on $\Xb_0$ which is not one of the $\tau_i$ would
correspond to a positive edge $e$ with $a_e\equiv 0\mod{m}$ such that
$\nu_e=1$ (by the assumption \eqref{eq3.1}).

To prove (iii), fix $i\in I$ and let us assume that $v:=t(e_i)$ is not
a leaf. Let $E_i$ be the set of edges with source $v$ which are
distinct from $\bar{e}_i$. There is a unique edge $e'\in E_i$ such
that $a_{e'}=a_i$, and $a_e\equiv 0\mod{m}$ for all $e\in E_i-\{e'\}$.
In other words, the $m$-cyclic cover $\Zb_v\to\Xb_v\cong\PP^1_k$ is
ramified at two points, so each component of $\Zb_v$ has genus $0$.
Since $\Yb_v\to\Zb_v$ is inseparable, each component of $\Yb_v$ has
genus $0$ as well, and $\Yb_v$ intersects with as many components of
$\Yb$ as $\Zb_v$ intersect with components of $\Zb$. Thus, the
minimality of the stable model implies that $|E_i|\geq 2$. Since
$\nu_{e_i}=\nu_i\in\{0,1\}$ (equation \eqref{eq3.4a}) we have
\begin{equation} \label{eq3.6}
  \sum_{e\in E_i}(\nu_e-1) \;=\;-1+\nu_i \in\{-1,0\}.
\end{equation}
It follows from (i) that $\nu_{e'}\geq 0$ and $\nu_e>1$ for all $e\in
E_i-\{e'\}$. We conclude that $E_i$ contains exactly two edges, $e'$
and $e''$. Furthermore, we find $\nu_i=1$, $\nu_{e'}=0$ and
$\nu_{e''}=2$. The remaining assertion that $t(e')=i$ and that
$t(e'')$ is a leaf is left to the reader (we will not use them in what
follows).  \Endproof

Here is an immediate consequence of Lemma \ref{lem3.1}:

\begin{cor} \label{cor3.1}
  Suppose that $f$ is not special. Then there exists an index $i\in I$
  such that the following holds. Let $\Xb_v$ be the component which
  meets the original component $\Xb_0$ in $\tau_i$. Then $\Xb_v$ has
  nontrivial intersection with exactly three components of $\Xb$
  (including $\Xb_0$). Furthermore, the torsor $Y\to Z$ has additive
  reduction over $\Xb_v$.
\end{cor}

See Figure \ref{treespic} for an illustration of the relevant part of
the tree $T$ and the curve $\Xb$, in the situation of Corollary
\ref{cor3.1}.

\subsection{The proof of Theorem \ref{thm2}}  \label{sec3.4}

In Section \ref{sec3.2} we have constructed, for any tuple of
$\Kb$-rational points $(x_i)$ lifting $(\tau_i)$, a metacyclic cover
$f:Y\to\PP^1$ of type $(x_i;a_i)$. In Section \ref{sec3.3} we have
analyzed the stable reduction of $f$. In this section we show that $f$
is special provided that $x_i\in\D_i$, with $\D_i$ as in
\eqref{diskeq}, thus proving Theorem \ref{thm2}. We do this in two
steps. First, we prove that $f$ is special if the $x_i$ are
$K_0$-rational. Note that this is a special case of Theorem
\ref{thm2}, as any $K_0$-rational point lifting $\tau_i$ is
automatically a center of the disk $\D_i$. Then, we show that the
cover $f$ remains special under a deformation which moves the branch
point $x_i$ into an arbitrary point $x_i'\in\D_i$.

\begin{prop}  \label{prop3.1}
  Let $x_1,\ldots,x_r$ be $K_0$-rational points of $\PP^1$ which lift
  $\tau_1,\ldots,\tau_r$. Then the metacyclic cover $f:Y\to\PP^1$ of
  type $(x_i;a_i)$ defined in Section \ref{sec3.2} is special.
\end{prop}

\proof By Proposition \ref{prop1.2}, the field of moduli of $f$ is
$K\inn=K_0(\zeta_p^{(m)})\subset K_0(\zeta_p)$. We may assume that the
field $K$ we have been working with in Sections \ref{sec3.2} and
\ref{sec3.3} is the minimal extension of $K\inn$ over which $f$ has
stable reduction. Recall from Section \ref{sec2.4} that the extension
$K/K\inn$ is Galois, of degree prime-to-$p$. The Galois group
$\Gamma:=\Gal(K/K\inn)$ acts faithfully and $k$-linearly on $\Yb$
(where $\fb:\Yb\to\Xb$ is the stable reduction of $f$), and this
action commutes with the action of $H$. We obtain an induced action of
$\Gamma$ on $\Zb$ and $\Xb$. The action of $\Gamma$ on the original
component $\Xb_0$ is trivial.

In our setup, the $p$-cyclic cover $Y\to Z$ reduces to a
$\bmu_p$-torsor $\Yb_0\to\Zb_0$ over the original component $\Xb_0$.
Therefore, the field $K$ contains the $p$th roots of unity.  Let
$\Gammat:=\Gal(K/K_0(\zeta_p))\subset\Gamma$. One shows easily that
$\Gammat$ acts trivially on $\Yb_0$ (see also the proof of Theorem
\ref{thm1b}).

Now suppose that $f$ is not special. Let $\Xb_v$ be a component as in
Corollary \ref{cor3.1}. In particular, $\Xb_v$ meets the rest 
of $\Xb$ in exactly three points. One of these points is $\tau_i$,
where $\Xb_v$ meets $\Xb_0$. Let $\tau'$ and $\tau''$ be the two other
points, corresponding to the edges $e'$ and $e''$ of Lemma
\ref{lem3.1} (iii). Note that $\tau_i$ and $\tau'$ are the two branch
points of the $m$-cyclic cover $\Zb_v\to\Xb_v$. Since the action of
$\Gammat$ is trivial on $\Xb_0$ and commutes with the $\ZZ/m$-action
on $\Zb$, it fixes all three points $\tau_i$, $\tau'$ and $\tau''$.
Therefore, $\Gammat$ acts trivially on $\Xb_v$.

Choose a point $\yb\in\Yb$ above $\tau_i$ and let $\zb$ be the image
of $\yb$ in $\Zb$. Note that $\yb$ and $\zb$ are ordinary double
points of $\Yb$ and $\Zb$, respectively. Let $\Yh_R$ be the completion
of $Y_R$ at $\yb$ and $\Zh_R$ the completion of $Z_R$ at $\zb$. Since
$Y_R$ is semistable, $\Yh_R$ is a {\em formal annulus}, of the form
$\Yh_R\cong\Spec R[[u,v\mid uv=\pi]]$, with $\pi\in R$. Corollary
\ref{cor3.1} states that we have additive reduction over the
component $\Zb_v$. By Remark \ref{rem2.1}, this implies 
\begin{equation} \label{eq3.8}
  0\;<\; v_R(\pi) \;<\; \frac{1}{(p-1)h_i}.
\end{equation}

Let $\lambda:=\zeta_p-1$. Note that $v_R(\lambda)=1/(p-1)$ and that
$\lambda$ is a uniformizer of $K_0(\zeta_p)$. We conclude from
\eqref{eq3.8} that $K_0(\zeta_p,\lambda^{1/N})\subset K$ and
$v_R(\pi)=c/(p-1)N$, with $N>h_i>m$ and $c$ prime-to-$N$.  Moreover,
$N$ is prime-to-$p$ because $K/K_0$ is of degree prime-to-$p$. Let
$\zeta_N$ be a primitive $N$th root of unity and $\sigma\in\Gammat$
such that $\sigma(\lambda^{1/N})=\zeta_N\lambda^{1/N}$. Then
$\sigma(\pi)=\zeta_N^c\pi$. Up to permutation, the parameter $u$
(resp.\ $v$) of the formal annulus $\Yh_R$ reduces to a local
parameter of $\Yb_0$ (resp. $\Yb_v$) at $\yb$.  Applying $\sigma$ to
the equation $uv=\pi$ and using the fact that $\Gammat$ acts trivially
on $\Yb_0$, we find that $\sigma$ induces an automorphism of $\Yb_v$
of order $N$ which fixes $\yb$. We conclude that $\sigma$ induces an
automorphism of $\Xb_v$ of order at least $N/(N,m)>1$. But this
contradicts the fact (proved earlier) that $\Gammat$ acts trivially on
$\Xb_v$. The proposition follows.  \Endproof

Let $f:Y\to\PP^1$ be as in Proposition \ref{prop3.1}. In particular,
$f$ is special. Proposition \ref{prop3.2} below states that, given
points $x_i'\in\D_i$, we can deform $f$ into a cover $f'$ of type
$(x_i';a_i)$. By construction, the cover $f'$ is special and gives
rise to $(\tau_i;a_i;\nu_i;\omega_0)$. This completes the proof of
Theorem \ref{thm2}.

\begin{prop} \label{prop3.2}
  For $i\in I$, choose $x_i'\in\D_i$. There exists a special cover
  $f':Y'\to\PP^1$ of type $(x_i';a_i)$ which gives rise to the special
  degeneration datum $(\tau_i;a_i;\nu_i;\omega_0)$. 
\end{prop}

\proof By Theorem \ref{thm1b}, the disk $\D_i$ (as defined by
\eqref{diskeq}, with center $x_i=\xt_i$) corresponds to the tail
$\Xb_i$, i.e.
\begin{equation} \label{eq3.9}
  \D_i \;=\; \{\; x\in\PP^1(\Kb) \mid\;
     \text{$x$ specializes to a point on 
                 $\Xb_i-\{\tau_i\}\cong\AA^1_k$}\;\}.
\end{equation}
In particular, $x_i'$ specializes to a point
$\xb_i'\in\Xb_i-\{\tau_i\}$. 

We define a finite, $H$-invariant map $\fb':\Yb'\to\Xb$ between
semistable $k$-curves. It is determined, up to unique isomorphism, by
the following requirements (compare with \cite{Raynaud98}, \S 3.2).
Over the original component $\Xb_0$, the maps $\fb:\Yb\to\Xb$ and
$\fb':\Yb'\to\Xb$ are the same. Even more, there exists an \'etale map
$U\to\Xb$ whose image contains $\Xb_0$ such that the pullbacks
$\fb|_U$ and $\fb'|_U$ are isomorphic, as $m$-cyclic covers of $U$.
Finally, for each $i\in I$, the restriction of $\fb'$ to the open
subset $\Xb_i-\{\tau_i\}\subset\Xb$ is a (possibly disconnected)
$H$-cover which is at most tamely ramified at $\xb_i'$.

Arguing as in \cite{Raynaud98}, \S 3.2, one shows that
$\fb':\Yb'\to\Xb$ lifts to a finite, $H$-invariant morphism
$f'_R:Y_R'\to X_R$ between semistable $R$-curves such that the
following holds. First, $f_R'$ is the stable model of an $H$-cover
$f':Y'\to\PP^1$ of type $(x_i';a_i)$. Second, the restrictions of
$f_R'$ and $f_R$ to the formal completion $\mathcal{X}_0$ of $X_R$
along $\Xb_0$ are isomorphic (as finite, $H$-invariant maps to
$\mathcal{X}_0$). It follows immediately that the $H$-cover $f'$ is
special and gives rise to the special degeneration datum
$(\tau_i;a_i;\nu_i;\omega_0)$. This proves Proposition \ref{prop3.2}
and therefore completes the proof of Theorem \ref{thm2}.
\Endproof

\subsection{The case $r=4$} \label{sec3.5}

Let us call two special degeneration data
$(\tau_i;a_i;\nu_i;\omega_0)$ and $(\tau'_i;a'_i;\nu'_i;\omega'_0)$
{\em equivalent} if there exists a $\ZZ/m$-equivariant isomorphism
$\phi:\Zb_0\iso \Zb_0'$ (where $\Zb_0\to\PP^1_k$ and
$\Zb_0'\to\PP^1_k$ are the corresponding $m$-cyclic covers) such that
$\phi^*\omega_0'=c\omega_0$, for some constant $c\in\FF_p^\times$. In
this section we determine all special degeneration data
$(\tau_i;a_i;\nu_i;\omega_0)$ with four branch points, up to
equivalence.

Fix $p$, $m$ and $(a_i)$, as in the beginning of Section \ref{sec1.1}.
It is clear that for given prime $p$, there is only a finite number of
possibilities for $m$ and $(a_i)$, which are easy to describe. After
reordering the indices, we may assume that $\nu_1=1$ and
$\nu_2=\nu_3=\nu_4=0$. Furthermore, after a projective linear
transformation, we may assume that $\tau_1=0$, $\tau_2=1$,
$\tau_4=\infty$. We write $\lambda$ instead of $\tau_3$; note that
$\lambda\in k-\{0,1\}$. Let $Z_\lambda\to\PP^1$ be the $m$-cyclic
cover given by the equation
\begin{equation} \label{eq3.5.1}
  z^m \;=\; x^{a_1}\,(x-1)^{a_2}\,(x-\lambda)^{a_3}.
\end{equation}
By an easy calculation, involving the divisors of $z$, $x$ and ${\rm
  d}x$, one shows that any differential $\omega_0$ on $Z_\lambda$
  satisfying Equation \eqref{eq3.1} is of the form
\begin{equation} \label{eq3.5.2}
  \omega_0 \;=\; \mu\;\frac{z\,{\rm d}x}{(x-1)(x-\lambda)}
           \;=\; \mu\;\frac{x\,{\rm d}x}{w},
\end{equation}
for some constant $\mu\in k^\times$. Here we have set
$w:=z^{-1}x(x-1)(x-\lambda)$. Note that 
$w^m=x^{a_1'}(x-1)^{a_2'}(x-\lambda)^{a_3'}$, with $a_i':=m-a_i$. 

Recall that the Cartier operator is defined
as the unique additive map $\C$ on differentials such that
\begin{align} \label{eqC1}
   \C(u^p\omega) &\;=\; u\,\C(\omega),\\ \label{eqC2}
   \C(\diff u)        &\;=\; 0,\;\; \text{\rm and} \\ \label{eqC3}
   \C(\diff u/u)      &\;=\; \diff u/u,
\end{align}
for all rational functions $u$ and differentials $\omega$ on
$Z_\lambda$.  We want to find all $\lambda$ and $\mu$ such that
$\C(\omega_0)=\omega_0$. We set $\alpha:=(p-1)/m$ and and
$f(x):=w^m=x^{a_1'}(x-1)^{a_2'}(x-\lambda)^{a_3'}$. Following
\cite{Yui80} and using \eqref{eqC1}--\eqref{eqC3}, we compute
\begin{equation} \label{eq3.5.3}
\begin{split}
   \C(\omega_0) &\;=\; \mu^{1/p}\;
   \C\big(\frac{x\,f(x)^\alpha\,\diff x}{w^p}\big)
                 \;=\;\mu^{1/p}\,w^{-1}\;
  \sum_{k=\alpha a_1'}^{\alpha(2m-a_1)}\,c_k^{1/p}\,\C(x^{k+1}\diff x)\\
                &\;=\; \mu^{1/p}\;(\,c_{p-2}^{1/p}\,
              +\,c_{2p-2}^{1/p}\,x\,)\,\frac{\diff x}{w},\\
\end{split}
\end{equation}
where $c_k$ is the coefficient of $x^k$ in $f(x)^\alpha$.  We see that
$\C(\omega_0)=\omega_0$ implies $c_{p-2}=0$. On the other hand, if
$c_{p-2}=0$, then $c_{2p-2}\neq 0$, because $\C(\omega_0)\neq 0$, by
Lemma \ref{lem1.0} (ii). Therefore, if $c_{p-2}=0$, we can set
$\mu:=c_{2p-2}^{1/(p-1)}$, to obtain $\C(\omega_0)=\omega_0$. We
conclude: if $\lambda$ is a zero of 
\begin{equation} \label{eq3.5.4}
  \Phi(\lambda) \;:=\; \sum_{l=0}^N\;
     \binom{\alpha a_2'}{N-l}\,\binom{\alpha a_3'}{l}\;\lambda^l
\end{equation}
(with $N:=p-\alpha a_1'-2=\alpha a_1-1$), then there exists a
differential $\omega_0$ on $Z_\lambda$ satisfying \eqref{eq3.1} and
\eqref{eq3.2}, unique up to a constant factor in $\FF_p^\times$. 

\begin{lem} \label{lem3.5.1}
  The zeros of the polynomial $\Phi$ in \eqref{eq3.5.4} are simple
  and $\neq 0,1$. 
\end{lem}

\proof It is shown in \cite{StvPvM} that the zeros of $\Phi$ other
than $0$ and $1$ are simple (the reason is that $\Phi$ satisfies a
certain hypergeometric differential equation). Hence we have to show
that $\Phi(0),\Phi(1)\neq 0$. Recall that $a_i'=m-a_i$ and $\sum_i
a_i=m$. Since $a_2'=a_1+2m-a_3'-a_4'\geq a_1$, we have $\alpha
a'_2>\alpha a_1-1=N$, so $\Phi(0)=\binom{\alpha a_2'}{N}\neq 0$. On
the other hand, we have
\begin{equation} \label{eq3.5.5}
  \Phi(1) \;=\; \text{\rm coeff.\ of $x^N$ in 
     $(x-1)^{\alpha(a'_2+a'_3)}$}\;=\; \binom{\alpha(a_2'+a_3')}{N}.
\end{equation}
But
$\alpha(a_2'+a_3')-N=\alpha(a_1'+a_2'+a_3'-m)+1=\alpha(m+a_4)+1>1$, so
$\Phi(1)\neq 0$ as well.
\Endproof

We can summarize the discussion as follows:

\begin{prop} \label{prop3.5.1}
  Let $r:=4$. Given $p$, $m$, $(a_i)$ and $(\nu_i)$ (satisfying the
  usual conditions, and with $\nu_1=1$), there are exactly $\alpha
  a_1-1$ nonequivalent special degeneration data
  $(\tau_i;a_i;\nu_i;\omega_0)$. In particular, for given $p$, there
  exist only a finite number of nonequivalent special degeneration
  data, and they are all defined over some finite field
  $\FF_{p^n}$. Therefore, we may assume $k=\bar{\FF}_p$ throughout.
\end{prop}

\begin{rem} \label{rem3.5.1}
  For $r >4$, it is still true that there exist only a finite number
  of equivalence classes of special degeneration data, for fixed $p$.
  This is less obvious than for $r=4$, because the polynomial $\Phi$
  is replaced by a system of $r-3$ equations in $r-3$
  variables. However, a deformation argument shows that this system
  has only finitely many solutions. It would be interesting to obtain
  a formula for the number of solutions, as in the case $r=4$.
\end{rem}

\bibliographystyle{abbrv} \bibliography{hurwitz}

\vspace{5mm}
\flushright{DRL, University of Pennsylvania\\
            209 South 33rd Street\\
            Philadelphia, PA 19104-6395\\
            wewers@math.upenn.edu}

\end{document}